\documentclass{amsart}
\usepackage{amsmath}
\usepackage{amssymb}
\usepackage{amsthm}
\usepackage{dsfont}
\usepackage{lineno}
\usepackage{subfigure}
\usepackage{graphicx}
\usepackage{multirow}
\usepackage{color}
\usepackage{xspace}
\usepackage{pdfsync}
\usepackage{hyperref} 
\usepackage{algorithmicx}
\usepackage{algpseudocode}
\usepackage{algorithm}
\graphicspath{{../SweepingFigures/}}

\newcommand{\bq}{\begin{equation}}
\newcommand{\eq}{\end{equation}}

\newcommand{\bO}{\mathcal{O}}

\newcommand{\Prop}{\mathcal{P}}

\newcommand{\Sm}{\mathcal{S}}

\algnewcommand{\LineComment}[1]{\State \(\triangleright\) #1}

\theoremstyle{lemma}

\theoremstyle{remark}

\newtheorem{remark}{Remark}

\begin{document}

\title[Fast Sweeping Methods]{Fast Sweeping Methods for Hyperbolic Systems of Conservation Laws at Steady State II}

\author{Bj\"orn Engquist}
\thanks{Bj\"orn Engquist, Department of Mathematics and ICES, The University of Texas at Austin, 1 University Station C1200, Austin, TX 78712
USA (engquist@math.utexas.edu).  
This author was partially supported by NSF DMS-1217203.
}

\author{Brittany D. Froese}
\thanks{Brittany D. Froese, Department of Mathematics and ICES, The University of Texas at Austin, 1 University Station C1200, Austin, TX 78712
USA (bfroese@math.utexas.edu).  
This author was partially supported by NSF DMS-1217203 and an NSERC PDF.
}

\author{Yen-Hsi Richard Tsai}
\thanks{Yen-Hsi Richard Tsai, Department of Mathematics and ICES, The University of Texas at Austin, 1 University Station C1200, Austin, TX 78712
USA (ytsai@math.utexas.edu).  
This author was partially supported by NSF DMS-1217203, NSF DMS-1318975, and a Simons Foundation Fellowship.  
}

\begin{abstract}
The idea of using fast sweeping methods for solving stationary systems of conservation laws has previously been proposed for efficiently computing solutions with sharp shocks.  We further develop these methods to allow for a more challenging class of problems including problems with sonic points, shocks originating in the interior of the domain, rarefaction waves, and two-dimensional systems.  We show that fast sweeping methods can produce higher-order accuracy.  Computational results validate the claims of accuracy, sharp shock curves, and  optimal computational efficiency.
\end{abstract}

\date{\today}

\maketitle

\section{Introduction}\label{sec:intro}

The numerical solution of systems of nonlinear (stationary) conservation laws,
\bq\label{eq:system}
\begin{cases}
 \nabla \cdot F(U) = a(U,x), & x\in\Omega\\
 B(U,x) = 0, & x\in\partial\Omega
\end{cases}
\eq
continues to be an important problem in numerical analysis.  A major challenge associated with this task is the need to compute non-classical solutions~\cite{CourantFriedrichs}, which leads to the need to develop numerical schemes that correctly resolve discontinuities in weak (entropy) solutions.  

\subsection{Related work}

A natural approach to computing these stationary solutions is to use an explicit time stepping or pseudo time stepping technique to evolve the system to steady state~\cite{AbgrallMezine_Steady,AbgrallRoe,ChouShuWENO,JiangShuWENO}.  
Several approaches are available for resolving shock fronts including front tracking schemes~\cite{Glimm_FrontTracking}, upstream-centered schemes for conservation laws (MUSCL)~\cite{Colella_MUSCL,vanLeer}, central schemes~\cite{JiangTadmor_CentralMultiD,NessyahuTadmor_Central}, essentially non-oscillatory (ENO) schemes~\cite{HartenENO}, and weighted essentially non-oscillatory (WENO) schemes~\cite{LiuOsherChan_WENO,OsherShu_ENO,Shu_ENO}. 
However, the computational efficiency of these schemes is restricted by a CFL condition and the need to evolve the system for a substantial time in order to reach the steady state solution.  In order to substantially improve the efficiency of these computations, it is desirable to develop methods that solve the steady state equations directly instead of through a time-evolution process.

Early work in this direction used Newton's method to solve a discrete version of the boundary value problem~\eqref{eq:system} using shock tracking techniques~\cite{GustafssonWahlund_SteadyFlow,ShubinStephens_SteadyShock}.
More recently, Newton solvers have been applied to WENO approximations of the steady Euler equations~\cite{HuLiTang_Newton}.  These Newton-based schemes result in faster methods, but in practice do not achieve linear computation time due to the low regularity of the discrete systems and the need for accurate initialisation.  For more general systems, a Gauss-Seidel scheme based on a Lax-Friedrichs discretisation of the steady state equations was described in~\cite{Chen_LFSweeping}. In~\cite{Shu_WENOHomotopy}, a homotopy approach was introduced to evolve from an initial condition to a steady state solution without the restriction of a CFL condition.

In~\cite{HongNozzle}, the authors characterised the structure of one-dimensional steady nozzle flow problems.  The analytical techniques they used involved combining smooth solution branches via the Rankine-Hugoniot conditions for a shock.  This work provided a better understanding of solutions of one-dimensional stationary systems of conservation laws, but did not provide a numerical method for numerically computing the solutions.

In~\cite{EFTSweeping}, we proposed a new fast sweeping approach for solving stationary conservation laws.  This approach was motivated by the fast sweeping methods that have been used for solving Hamilton-Jacobi equations~\cite{KaoOsherTsai_Sweeping,TsaiChenOshwerZhao_Sweeping,ZhaoSweeping} and involves combining smooth solution branches as in~\cite{HongNozzle}.  These sweeping methods exploit the flow of information along characteristics, which allows solutions to be computed by passing through the domain in a small number of pre-determined sweeping directions.  The resulting algorithms typically have linear computational complexity.  

In our earlier work we proposed  a basic two-step sweeping method for systems of conservation laws.
\begin{enumerate}
\item Smooth solution branches are generated by means of an update formula that is used to update the solution along different sweeping directions.  
\item A selection principle based on the Rankine-Hugoniot and entropy conditions for a stationary shock is used to determine which solution branch is active at each point.  
\end{enumerate}
We used this framework to compute the solution to several simple one- and two-dimensional problems.  In two dimensions, the appropriate implementation of these two steps was informed by some \emph{a priori} knowledge of the basic structure of the solutions.  The resulting methods also had optimal computational complexity and produced sharp shocks.

\subsection{Contributions of this article}
In this article, we develop the basic idea of fast sweeping methods in order to construct methods for solving a larger range of problems in one- and two-dimensions.  These new developments allow us to demonstrate additional advantages of sweeping over more traditional shock-capturing methods.  

In \autoref{sec:accuracy}, we show how to construct higher order sweeping methods, which preserve their global accuracy even in the presence of shocks and sonic points.  This is a clear advantage over formally higher-order shock-capturing methods, which can be restricted to first-order accuracy when shock are present~\cite{EngquistSjogreen}.

In \autoref{sec:matching}, we consider the problem of sweeping and matching techniques in two dimensions.  In our earlier work, we proposed a technique for sweeping and matching when the exact solution is composed of smooth solution branches separated by a shock curve.  We now develop these techniques to allow us to solve problems with a more complicated characteristic structure.  In one of the problems we study, a shock begins in the middle of the domain rather than at the boundary. A second example involves rarefaction waves.  These methods produce sharp shock curves and edges, in distinction to the smeared out results produced by capturing methods.  In a second example, we demonstrate the use of sweeping methods for problems with rarefaction waves.

Finally, in \autoref{sec:systems}, we focus on the solution of a class of two-dimensional systems.  In this setting, we cannot view sweeping and matching as distinct steps since only partial boundary data may be available in some regions.  To solve these systems, we combine the usual selection principle (Rankine-Hugoniot conditions) with a condition that solutions are regular on either side of the shock.  In this section, we show how to construct an efficient method that simultaneously constructs a smooth solution state and a clean shock curve.  

\section{Accuracy}\label{sec:accuracy}
The development of higher-order schemes for solving systems of conservation laws has been an important area of research, with the essentially non-oscillatory (ENO)~\cite{HartenENO,ShuOsher_ENO} and weighted essentially non-oscillatory (WENO)~\cite{LiuOsherChan_WENO} schemes becoming very popular approaches.  When solutions include shocks, these methods will naturally produce an $\bO(1)$ error near the shock.  Depending on the structure of characteristics in the system, this error can propagate and cause the global solution accuracy to be only first-order even in regions where the solution is smooth~\cite{EngquistSjogreen}.

The fast sweeping methods involve computing smooth solution branches, which are joined together by directly imposing the Rankine-Hugoniot conditions at a shock.  The methods are also very flexible in the sense that we can solve the necessary ODEs or PDEs using any consistent, stable scheme.  This allows us to easily compute both smooth solution branches and the shock location with higher accuracy.  In particular, we can compute higher-order solutions in settings where traditional shock-capturing methods are limited to first-order accuracy.

\subsection{Solutions with Shocks}
In this section we describe a higher-order left-to-right sweeping method; the opposite sweeping direction can be handled in the same way.  

In one-dimensions, systems take the form
\bq\label{eq:CL1d}
\begin{cases}
f(U)_x = a(U,x), & x_L < x < x_R\\
U(x) = U_L, & x = x_L\\
B_R(U(x)) = 0, & x = x_R.
\end{cases}
\eq
We introduce the operator $\Prop_{x_Lx}U_L$ to denote the solution $U$ at the point $x$.  This is obtained by propagating the left boundary condition $U_L$ from $x_L$ to $x$ via the solution of the ODE
\bq\label{eq:ODE1d}
\begin{cases}
V_x = a(f^{-1}(V),x), & x > x_L\\
V(x) = f(U_L), & x = x_L.
\end{cases}
\eq
We also introduce the jump operator $\Phi U_-$, which gives an entropy-satisfying solution of the Rankine-Hugoniot conditions
\bq\label{eq:jump}
f(\Phi U_-) = f(U_-).
\eq
When the solution consists of two smooth states separated by a single shock, the sweeping method requires solving the following equation for the unknown shock location $x_S$:
\bq\label{eq:sweep1D}
B_R(\Prop_{x_Sx_R}^h\Phi\Prop_{x_Lx_S}^hU_L) = 0.
\eq
Here $\Prop^h$ refers to a discrete version of the propagation operator and $h$ denotes the stepsize on the grid.

Theorem~1 in~\cite{EFTSweeping} ensures that if the method used to solve the necessary ODEs is consistent, stable, and has accuracy $\bO(h^k)$, then the solution computed by the sweeping method will also have accuracy $\bO(h^k)$ in $L^1$.  As the first example in this paper, we demonstrate that higher-order ODE solvers can easily be incorporated into the method.  Here we provide an example where traditional shock-capturing methods are only first-order accurate, while the fast sweeping method successfully achieves higher accuracy.

\subsubsection{Nozzle problem (shock)}
We consider a nozzle problem that was discussed in~\cite{EngquistSjogreen}. In that work, the authors solve this problem using a fourth-order accurate ENO method and a second-order accurate TVD scheme.  They consider the $L^2$ error in a region to the left ($0<x<4.5$) and right ($5.5 < x < 10$) of the shock.  In the left region, solutions obtain the expected accuracy.  However, to the right of the shock, the accuracy is only first-order; this is a consequence of the characteristic structure at the shock.

The system we want to solve is
\bq\label{eq:nozzle} 
\left(\begin{tabular}{c}$\rho A$\\ $\rho u A$\\ $EA$\end{tabular}\right)_t + \left(\begin{tabular}{c}$\rho u A$\\ $(\rho u^2 + p) A$\\ $uA(E+p)$\end{tabular}\right)_x = \left(\begin{tabular}{c}0\\$pA'(x)$\\ 0\end{tabular}\right), 
\eq
which we want to solve to steady state on the domain $x\in[0,10]$.  Here the nozzle area is
\[ A(x) = 1.398 + 0.347 \tanh(0.8x-4), \]
the pressure is
\[ p = (\gamma-1)\left(E-\frac{1}{2}\rho u^2\right) = \rho R T, \]
and we take $\gamma=1.4$ and $R = 8.3144$.
We enforce the boundary conditions
\[ \rho=0.502, \quad u = 1.299, \quad p = 0.3809 \]
at $x=0$ and
\[ \rho = 0.7519 \]
at $x=10$.

\begin{figure}
\centering
\includegraphics[width=0.6\textwidth]{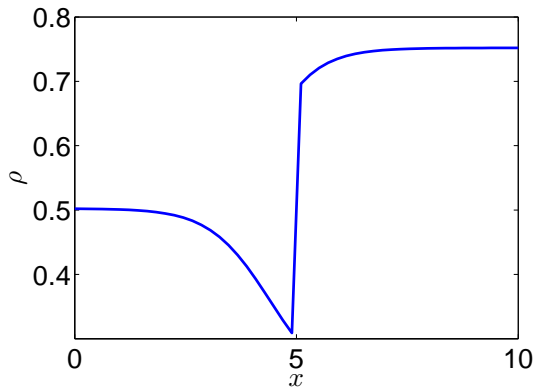}
\caption{Density at steady state for nozzle problem.}
\label{fig:densityNozzleShock}
\end{figure}

We can integrate the system to obtain an implicit expression for the exact solution, which can be approximated to machine precision.  This gives us an exact solution to use to compute the error in our approximations.  We note that the exact solution contains a shock at $x\approx5$ (Figure~\ref{fig:densityNozzleShock}).

\begin{figure}
\centering
\subfigure[]{\includegraphics[width=0.8\textwidth]{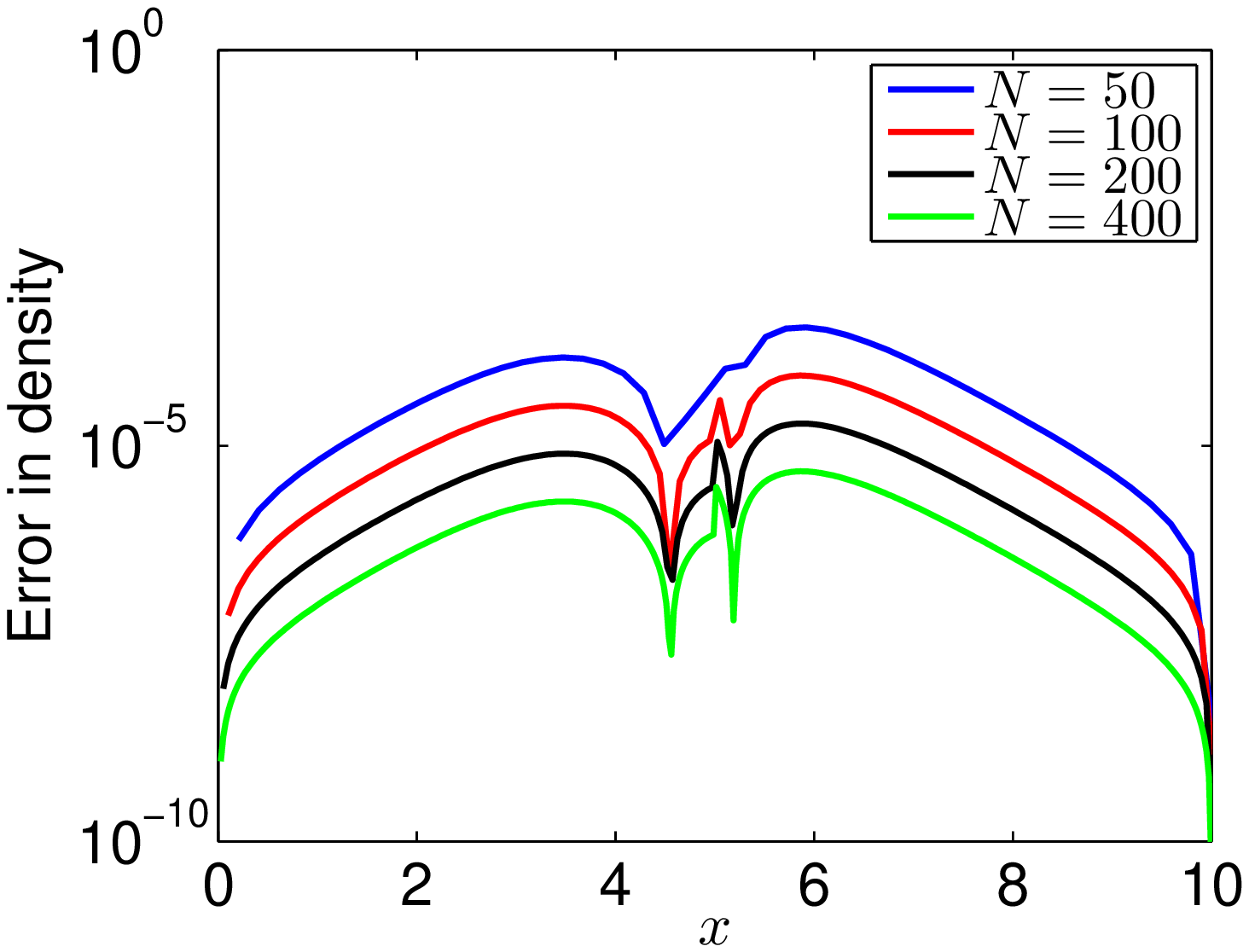}\label{fig:errorTrap}}
\subfigure[]{\includegraphics[width=0.8\textwidth]{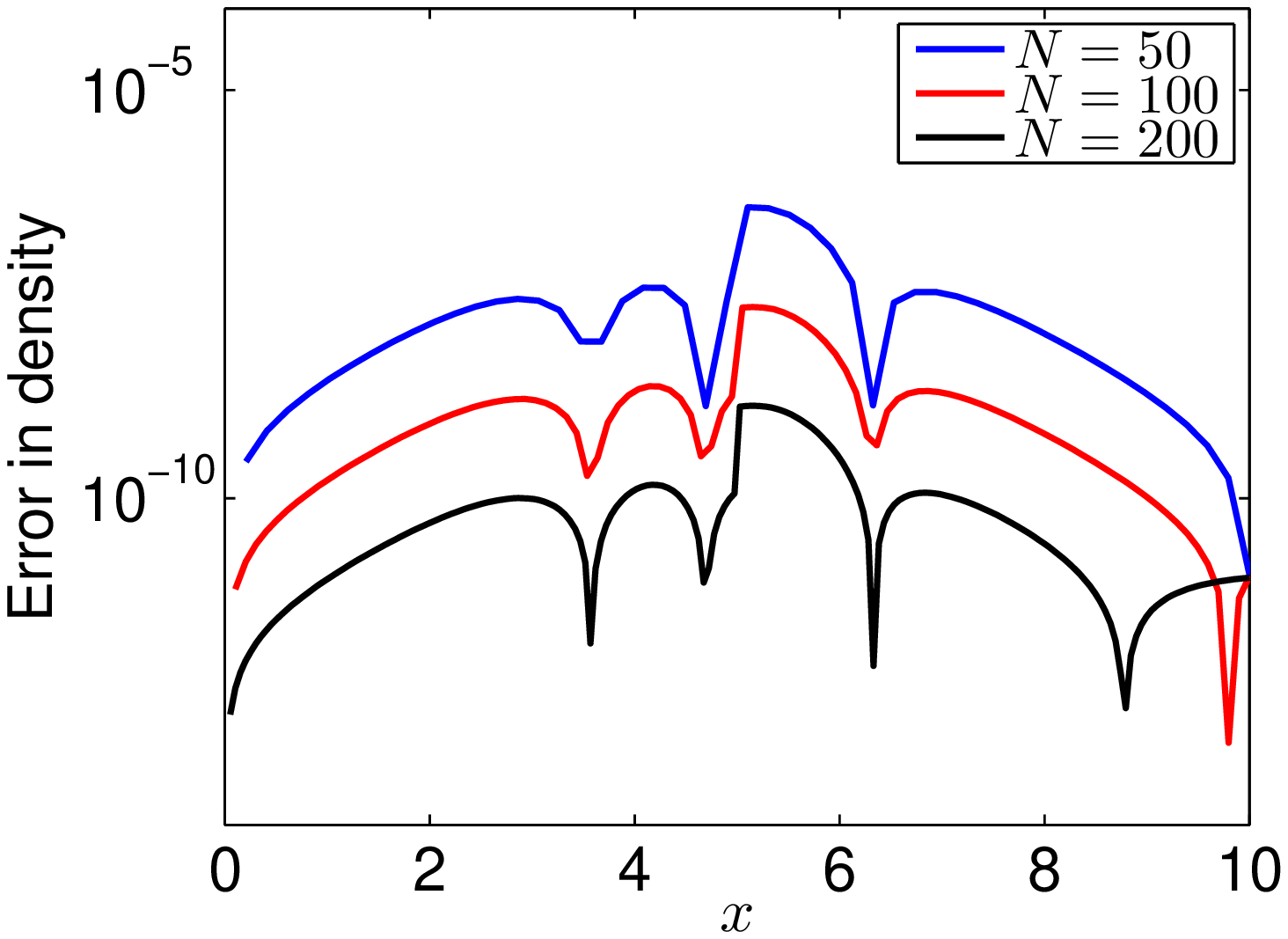}\label{fig:errorRK4}}
\caption{Error in density for the solution of the nozzle problem by sweeping with \subref{fig:errorTrap}~the trapezoid rule and \subref{fig:errorRK4}~RK4.}
\label{fig:errorDensity}
\end{figure}

We perform the sweeping (integration) using two different methods: the trapezoid rule and a fourth-order Runge-Kutta method.  We plot the error in the density in Figure~\ref{fig:errorDensity}. We also provide the $L^2$ solution error in the intervals $[0,4.5]$ and $[5.5,10]$ (to the left and right of the shock); see Table~\ref{table:errorTrap}.  In distinction to shock-capturing methods, we find that the fast sweeping methods preserve the higher-order accuracy on both sides of the shock.

\begin{table}[htdp]
\caption{$L^2$ error in density to the left and right of the shock for the solution of the nozzle problem by sweeping.}
\begin{center}
\begin{tabular}{cc|cc|cc|cc}
&&\multicolumn{2}{c}{$0<x<4.5$}&\multicolumn{2}{c}{$5.5<x<10$}&\multicolumn{2}{c}{Shock Location}\\
Scheme & $N$  & Error in $L^2$ & Rate & Error in $L^2$ & Rate & Error & Rate\\
\hline
Trap & 50  & $6.663\times10^{-5}$ &   &$1.438\times10^{-4}$&   & $1.412\times10^{-3}$ &\\
Trap & 100 & $1.149\times10^{-5}$ &2.5&$0.242\times10^{-4}$&2.6& $0.428\times10^{-3}$ &1.7\\
Trap & 200 & $0.201\times10^{-5}$ &2.5&$0.042\times10^{-4}$&2.5& $0.116\times10^{-3}$ &1.9\\
Trap & 400 & $0.035\times10^{-5}$ &2.5&$0.007\times10^{-4}$&2.5& $0.030\times10^{-3}$ &1.9\\
\hline
RK4 & 50  & $1.829\times10^{-8}$ &   &$7.934\times10^{-8}$&   & $2.017\times10^{-6}$ &\\
RK4 & 100 & $0.076\times10^{-8}$ &4.6&$0.276\times10^{-8}$&4.8& $0.113\times10^{-6}$ &4.2\\
RK4 & 200 & $0.003\times10^{-8}$ &4.5&$0.012\times10^{-8}$&4.5& $0.007\times10^{-6}$ &4.1\\
\end{tabular}
\end{center}
\label{table:errorTrap}
\end{table}

\subsection{Solutions with turning points}
A more challenging situation occurs when the solution contains a turning points; that is, when one (or more) of the eigenvalues of the flux gradient ($\nabla f$) changes sign smoothly.  Here we show how to develop higher-order fast sweeping methods for solving one-dimensional systems that contain turning points.

We describe in detail the solution of a problem with a single turning point, where an eigenvalue smoothly transitions from negative to positive.  The same technique can be used to solve problems with multiple turning points, as outlined in~\cite{EFTSweeping}.  

Note that this approach is only needed when an eigenvalue transitions from negative to positive through the turning point.  In the  situation where an eigenvalue transitions continuously from positive to negative, the turning point can be interpreted as a degenerate shock; the jump operator will have a continuous solution that satisfies the Lax entropy conditions~\cite{Lax}.  In this case we can simply use the Rankine-Hugoniot conditions and entropy conditions to determine the turning point location.

Recall that at a turning point $x_T$, where an eigenvalue $\lambda_i$ vanishes, the solution $U_T$ must satisfy a compatibility condition
\bq\label{eq:solvability}  (P^{-1}a\left(U_T,x_T\right))_i = 0. \eq
Here the matrix $P$ is defined via a decomposition of the Jacobian of the flux,
\bq\label{eq:defP}
P^{-1}\nabla f(U_T)P = \Lambda,
\eq
where $\Lambda$ is a diagonal matrix containing the eigenvalues of $\nabla f$.  For a well-posed problem, the presence of a turning point will be accompanied by ``incomplete'' boundary conditions.  That is, instead of providing values for all solution components, the given boundary condition on the left will be a curve through state space that is parameterised by an unknown $\alpha$,
\[ U(x_L) = U_L^\alpha. \]
Then the unknown boundary condition is determined by compatibility condition at the turning point.  That is, we search for the unknown $\alpha$ that solves the equation
\bq\label{eq:unknownBC}  (P^{-1}a(\Prop_{x_Lx_T}U_L^\alpha,x_T))_i = 0.\eq

In these situations, we can still use a conventional ODE solver to solve the system of ODES from $x_L$ up to a grid point $x_j<x_T$ that is near the sonic point.

Next, we need to solve the following system for the unknown solution $U_T$ and the turning point location $x_T = x_j + h$.
\bq\label{eq:turn_HO}
\begin{cases}
U_T=\Prop^h_{x_j,x_T}U_j\\
\lambda_i(U_T) = 0.
\end{cases}
\eq
In the above equation, we require that $\lambda_i$ is the largest negative eigenvalue since the eigenvalue is smoothly transitioning from negative to positive.

In~\cite{EFTSweeping}, we used a forward Euler step to approximate both the solution at the turning point and the location of the turning point. However, this approach does not generalise naturally to higher-order methods because of the (possible) need to accurately invert $\nabla f$ at or near a point where an eigenvalue vanishes. 
Instead, we will rely on the regularity of the solution near the turning point.

We recall that we have the values of both $U$ and $U_x = (\nabla f)^{-1}a$ at the grid point $x_k, k\leq j$.  We will use these values to compute $U_{j+1}$ by extrapolation.  For example, if we let $\Delta x$ = $x_j-x_{j-1}$ and $0 \leq h \leq \Delta x$, we can use the formula
\bq\label{eq:extrap}
\begin{split}
U(x_j + h) = &\left(1+\frac{h}{\Delta x}\right)U(x_j) + \frac{1}{2}h\left(1+\frac{h}{\Delta x}\right)U_x(x_j)\\ &- \frac{h}{\Delta x}U(x_{j-1}) - \frac{1}{2}h\left(1+\frac{h}{\Delta x}\right)U_x(x_{j-1}) + \bO(\Delta x^4).
\end{split}
\eq
Similarly, we can use more or fewer of the neighbouring points in order to obtain an approximation that has the desired order of accuracy.  Notice that this formula can be interpreted as a non-traditional linear multistep method.

We will use this formula for two purposes:
\begin{itemize}
\item It is used to compute $U_{j+1} $without the need for inverting the flux function.
\item It is used to determine the location of the turning point.  That is, we find the value of $h$ for which $\lambda(U(x_j + h)) = 0.$  The turning point is then given by $x_T = x_j + h$.
\end{itemize}

This gives us an approach for propagating a solution from the left boundary $x_L$ all the way up to the turning points $x_T$.  Using this method, we can solve~\eqref{eq:unknownBC} for the ``missing'' boundary conditions.  From here, we can continue to use a conventional ODE solver to propagate the solution from $x_{j+1}$ to the right end of the domain.

\subsubsection{Nozzle problem (sonic point and shock)}
We provide some computational results for the nozzle problem that we computed in~\cite{EFTSweeping} using forward Euler.  
\bq\label{eq:nozzle} 
\left(\begin{tabular}{c}$\rho A$\\ $\rho u A$\\ $EA$\end{tabular}\right)_t + \left(\begin{tabular}{c}$\rho u A$\\ $(\rho u^2 + p) A$\\ $uA(E+p)$\end{tabular}\right)_x = \left(\begin{tabular}{c}0\\$pA'(x)$\\ 0\end{tabular}\right), 
\eq
which we want to solve to steady state on the domain $x\in[0,3]$.

Here the pressure is
\[ p = (\gamma-1)\left(E-\frac{1}{2}\rho u^2\right) = \rho R T \]
and the sound speed is
\[ c = \sqrt{\gamma p/\rho}. \]
The eigenvalues of the Jacobian are $\lambda_1 = u-c$, $\lambda_2 = u$, and $\lambda_3 = u+c$.

We take the cross-sectional area to be
\[ A(x) = 1+2.2(x-1.5)^2, \]
the gas constant $\gamma = 1.4$, and $R = 8.3144$.
We consider the boundary conditions 
\[ p_L = 1, \quad p_R = 0.6784, \quad T_L  = 300. \]
Boundary data for the remaining solution components is computed via compatibility conditions as detailed in~\cite[section 3.5.1]{EFTSweeping}.

\begin{figure}[htd]
	\centering
			\subfigure[]{\includegraphics[width=.48\textwidth]{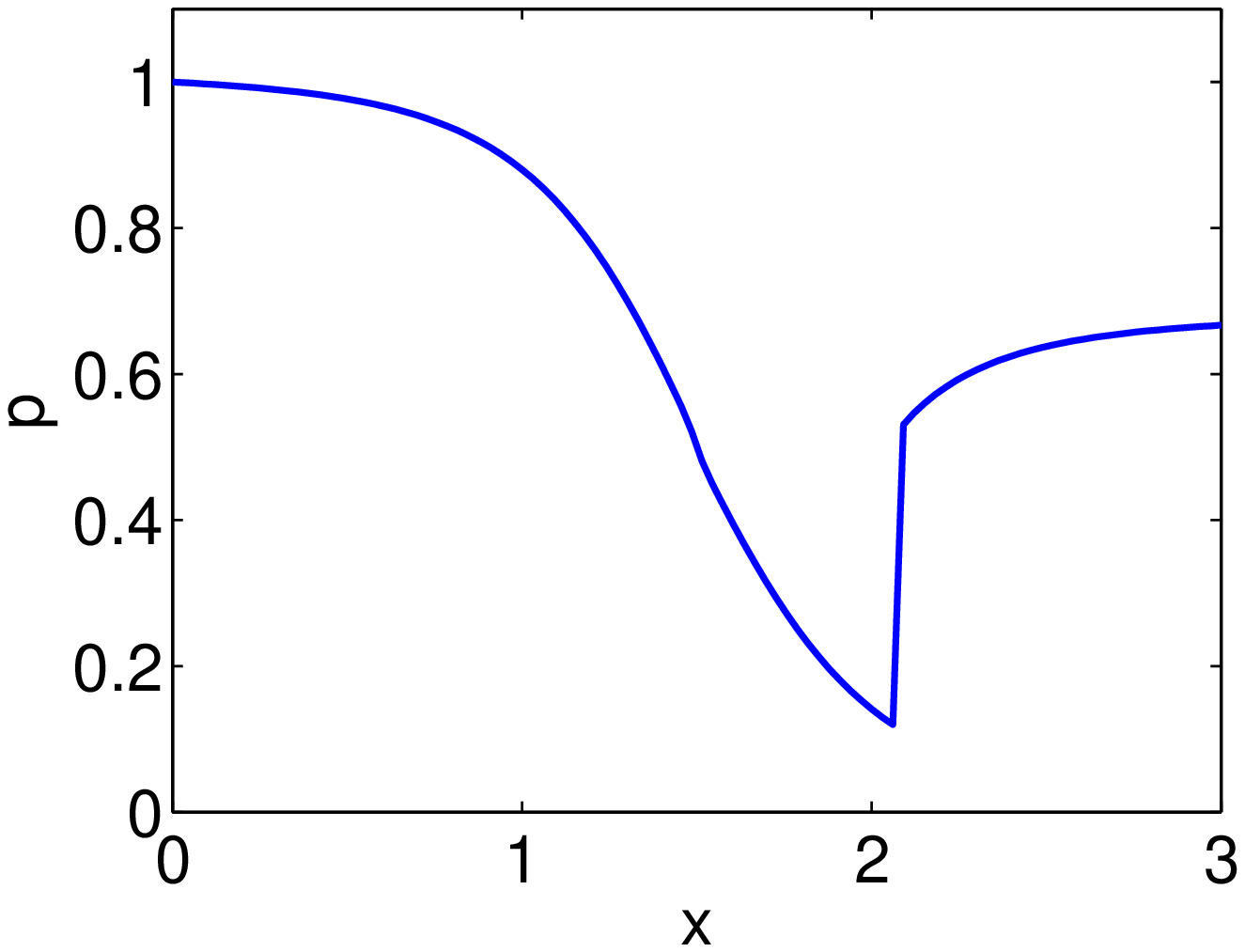}\label{fig:p_nozzle_100}}
       \subfigure[]{\includegraphics[width=.48\textwidth]{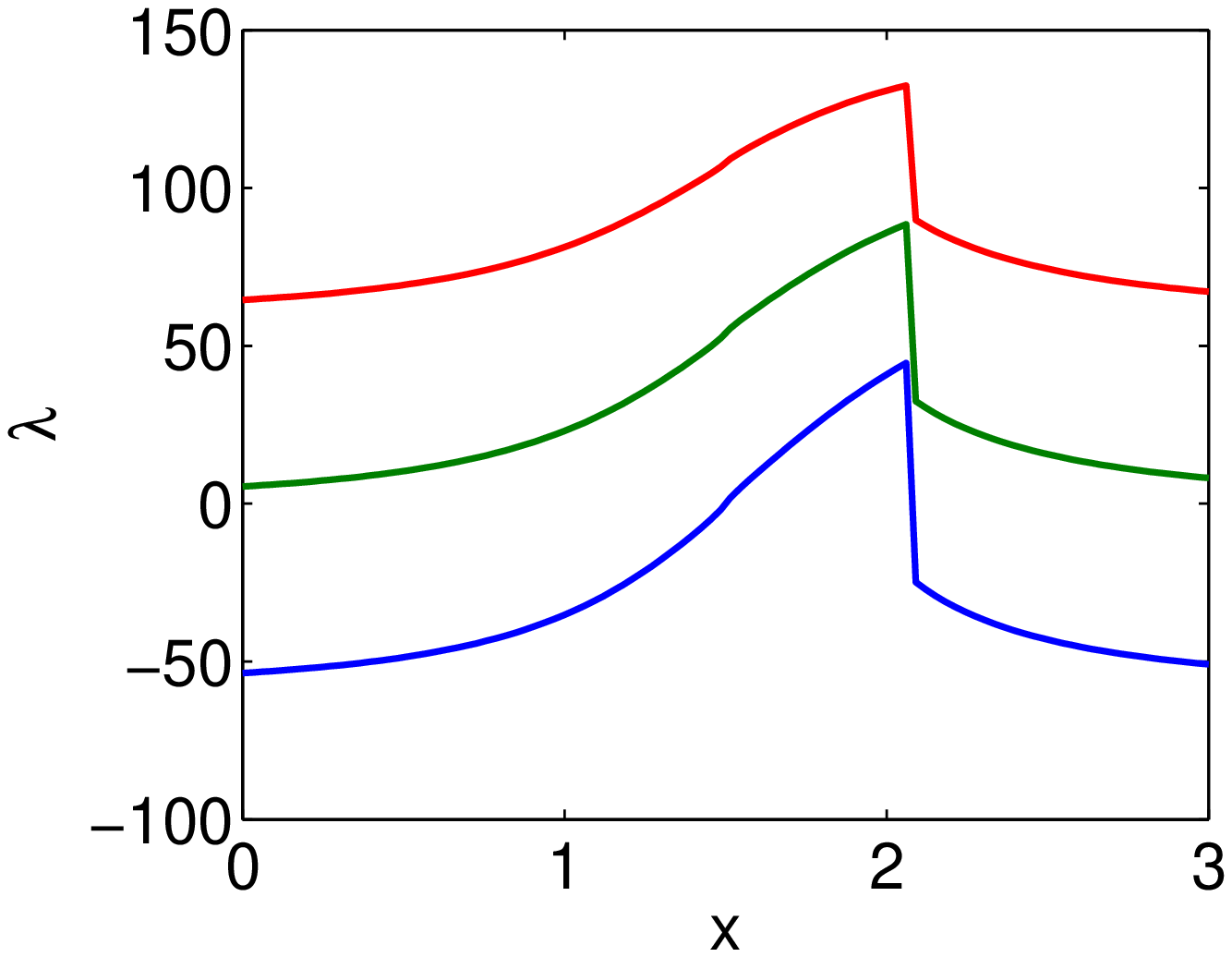}\label{fig:l_nozzle_100}}
  	\caption{Solution to the nozzle problem: \subref{fig:p_nozzle_100}~pressure and \subref{fig:l_nozzle_100}~eigenvalues.}
  	\label{fig:nozzle}  	
\end{figure}

\begin{figure}
\centering
\subfigure[]{\includegraphics[width=0.8\textwidth]{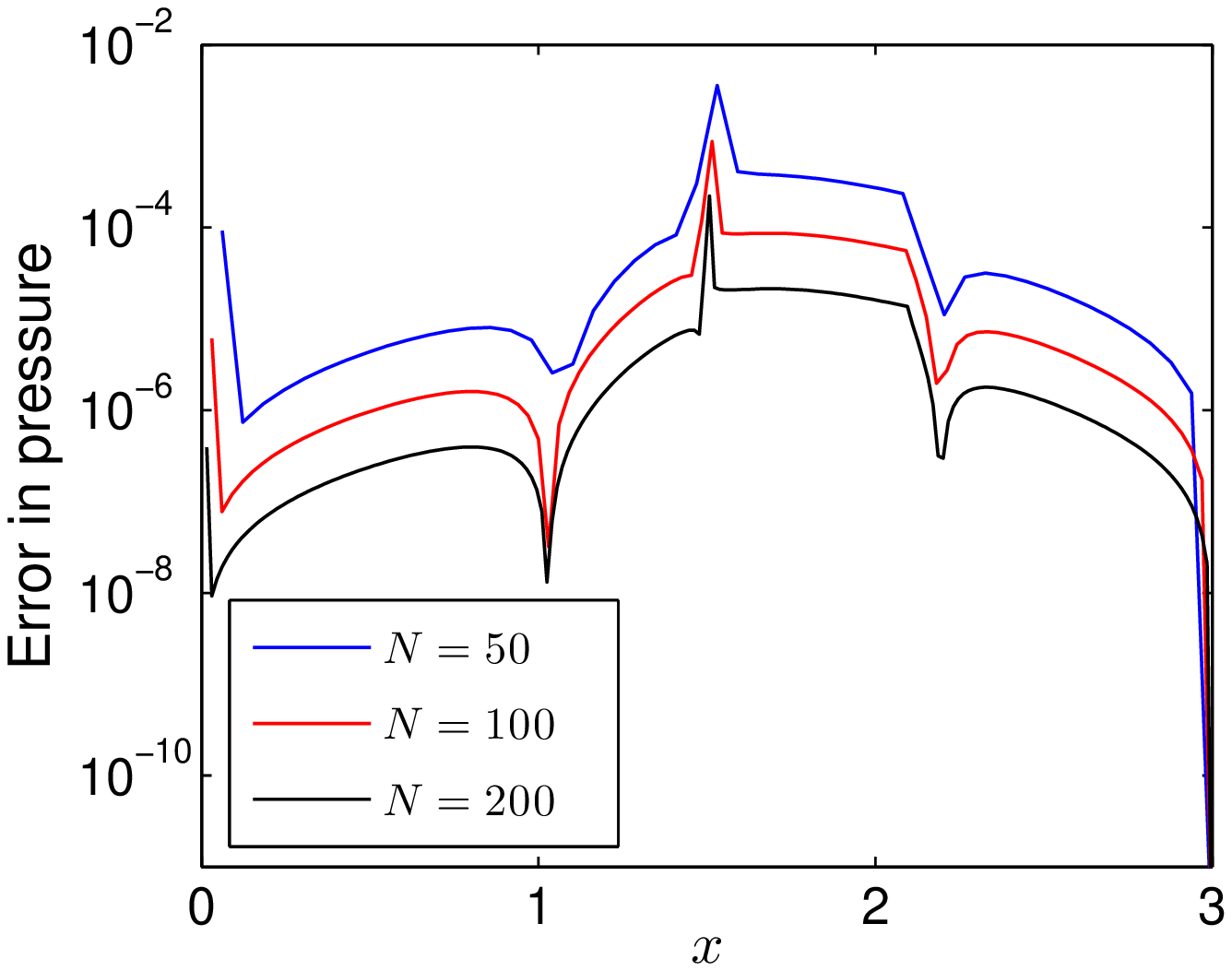}\label{fig:errorSonTrap}}
\subfigure[]{\includegraphics[width=0.8\textwidth]{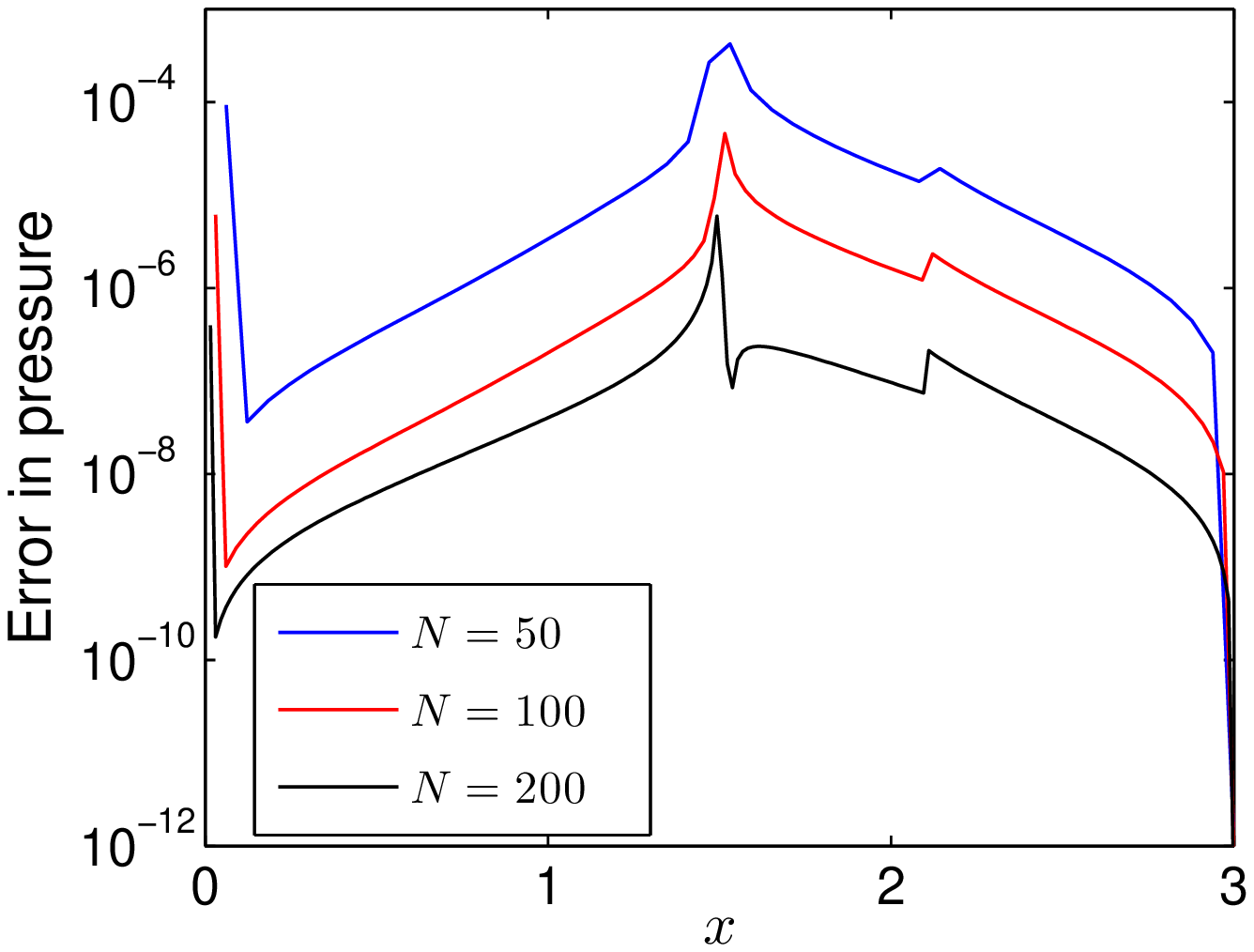}\label{fig:errorSonRK4}}
\caption{Error in pressure for the solution of the nozzle problem by sweeping with \subref{fig:errorSonTrap}~the trapezoid rule and \subref{fig:errorSonRK4}~RK4.}
\label{fig:errorSon}
\end{figure}

The solution (pressure and eigenvalues) are in Figure~\ref{fig:nozzle}.  We compute the solution using the trapezoid rule and a fourth-order Runge-Kutta method.  The error is presented in Figure~\ref{fig:errorSon} and Table~\ref{table:errorSon}.  Despite the challenge of this problem, we successfully compute the solution to the desired accuracy.

\begin{table}[htdp]
\caption{$L^2$ error in pressure for the solution by sweeping of the nozzle problem with a sonic point and shock.}
\begin{center}
\begin{tabular}{cc|cc}
Scheme & $N$  & Error in $L^2$ & Rate \\
\hline
Trap & 50  & $2.292\times10^{-4}$ &   \\
Trap & 100 & $0.287\times10^{-4}$ &2.9\\
Trap & 200 & $0.039\times10^{-4}$ &2.7\\
\hline
RK4 & 50  & $3.318\times10^{-5}$ &   \\
RK4 & 100 & $0.164\times10^{-5}$ &4.3\\
RK4 & 200 & $0.010\times10^{-5}$ &4.0\\
\end{tabular}
\end{center}
\label{table:errorSon}
\end{table}

\section{Shock Capturing versus Matching}\label{sec:matching}
In~\cite{EFTSweeping}, we also outlined a basic sweeping approach for two-dimensional problems, which was tested on several simple examples.  The approach consisted of two basic steps:
\begin{enumerate}
\item Compute smooth solution branches by sweeping in the boundary conditions.
\item Use the Rankine-Hugoniot conditions to generate a shock curve that separates two smooth states.
\end{enumerate}

In some situations, it is possible to generate a solution (including shocks) by sweeping in one direction only.  This leads to a simple method for solving stationary systems that can have a fairly complicated characteristic structure.  However, it also negates one of the benefits of fast sweeping methods, which is the ability to produce sharp shock curves by matching smooth states.  

Now we show how the basic sweeping method can be used to construct solutions that have more complicated charasteristic structures.  This essentially involves breaking the problem into small sub-domains, where a simple fast sweeping method can be applied.

\subsection{Shock curve generated in interior of domain}

One example where sweeping is desirable but non-trivial is the example
\bq\label{eq:2dscalar} \left(\frac{u^2}{2}\right)_x + u_y = 0 \eq
subject to the boundary conditions 
\[ u(0,y) = 1.5, \quad u(1,y) = -0.5, \quad u(x,0) = 1.5-2x. \]
In this example, $f(u)=\frac{1}{2}u^2$ and $g(u) = u$.

\begin{figure}[htdp]
\begin{center}
\subfigure[]{\includegraphics[width=0.45\textwidth]{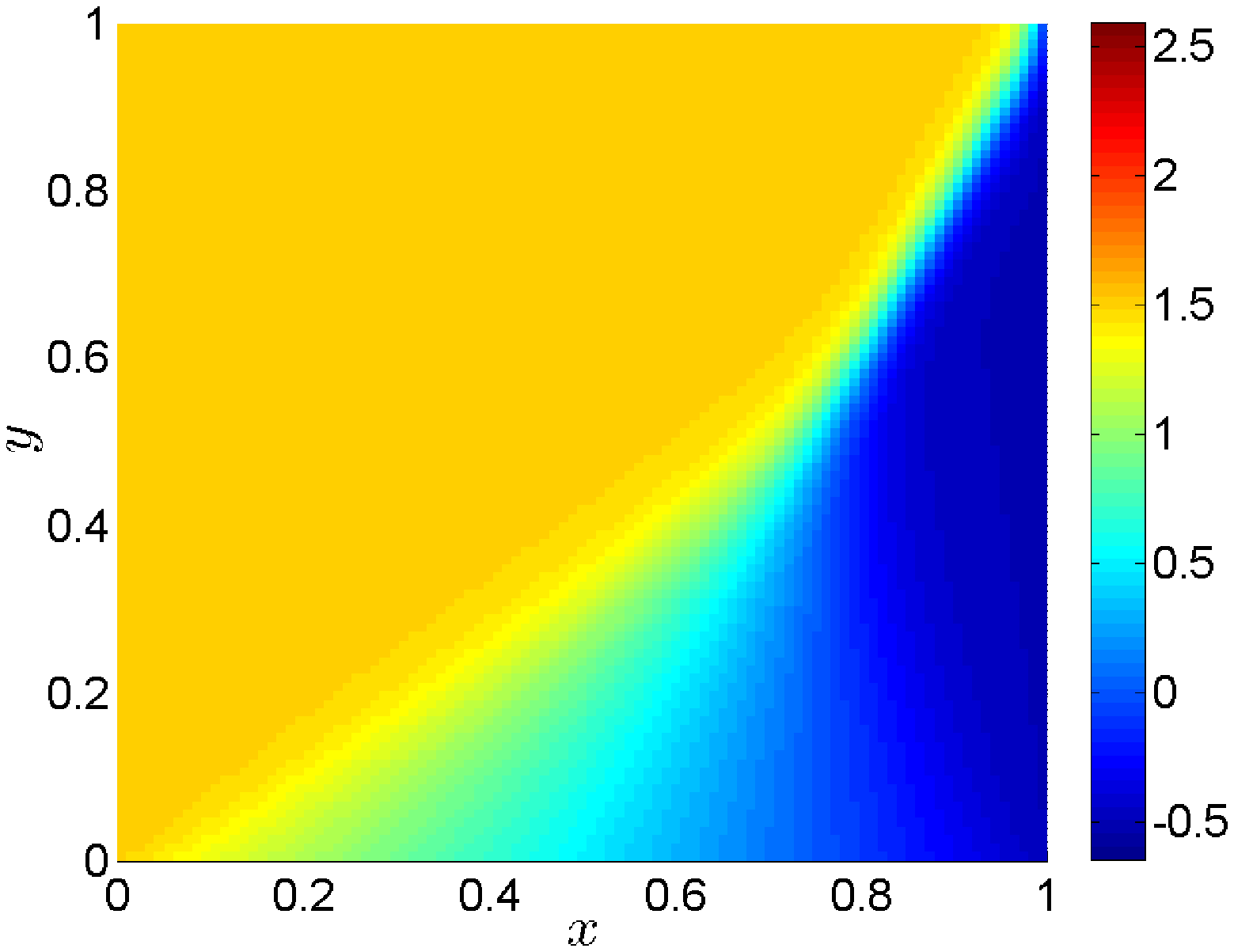}\label{fig:stupidLF}}
\subfigure[]{\includegraphics[width=0.45\textwidth]{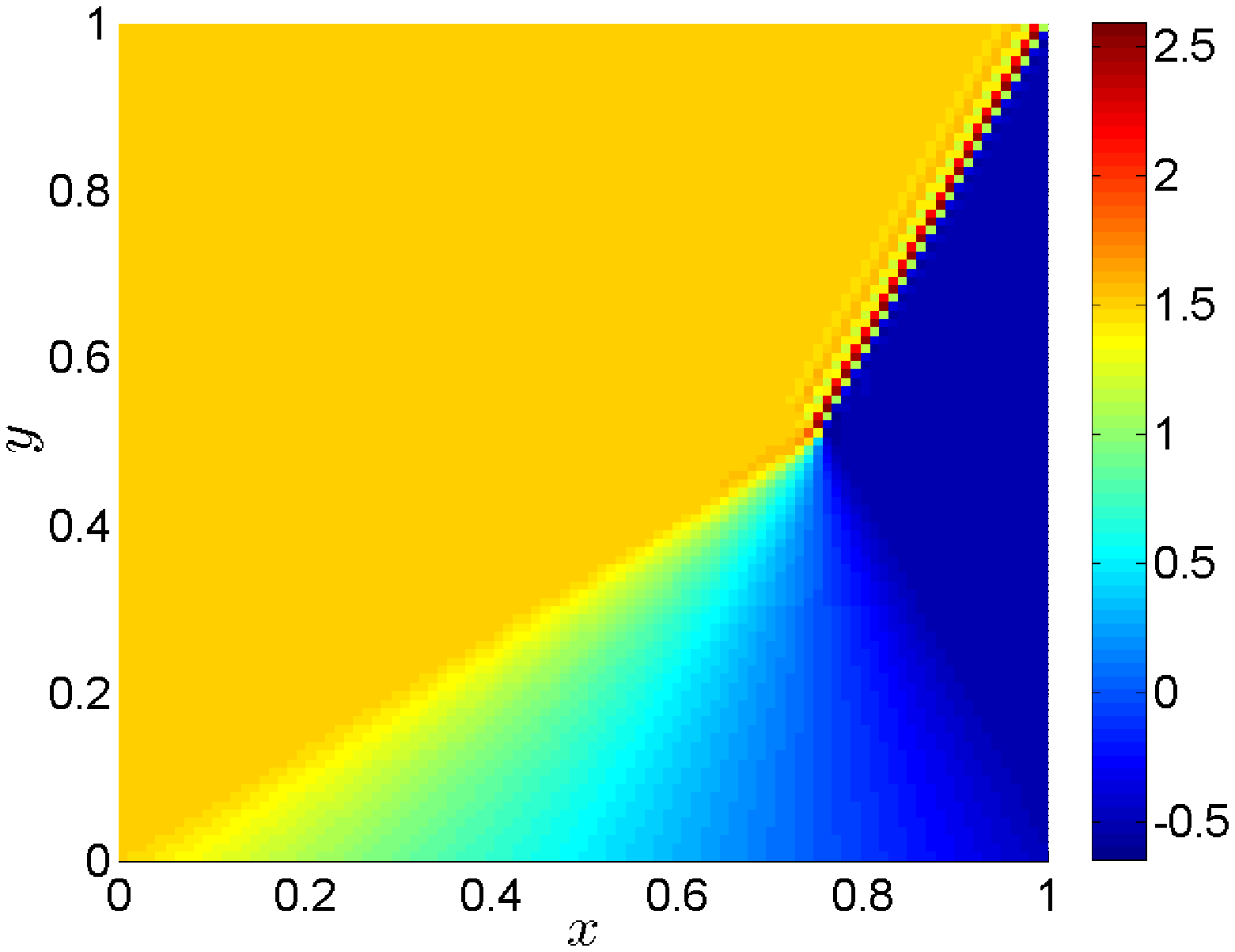}\label{fig:stupidLW}}\\
\subfigure[]{\includegraphics[width=0.45\textwidth]{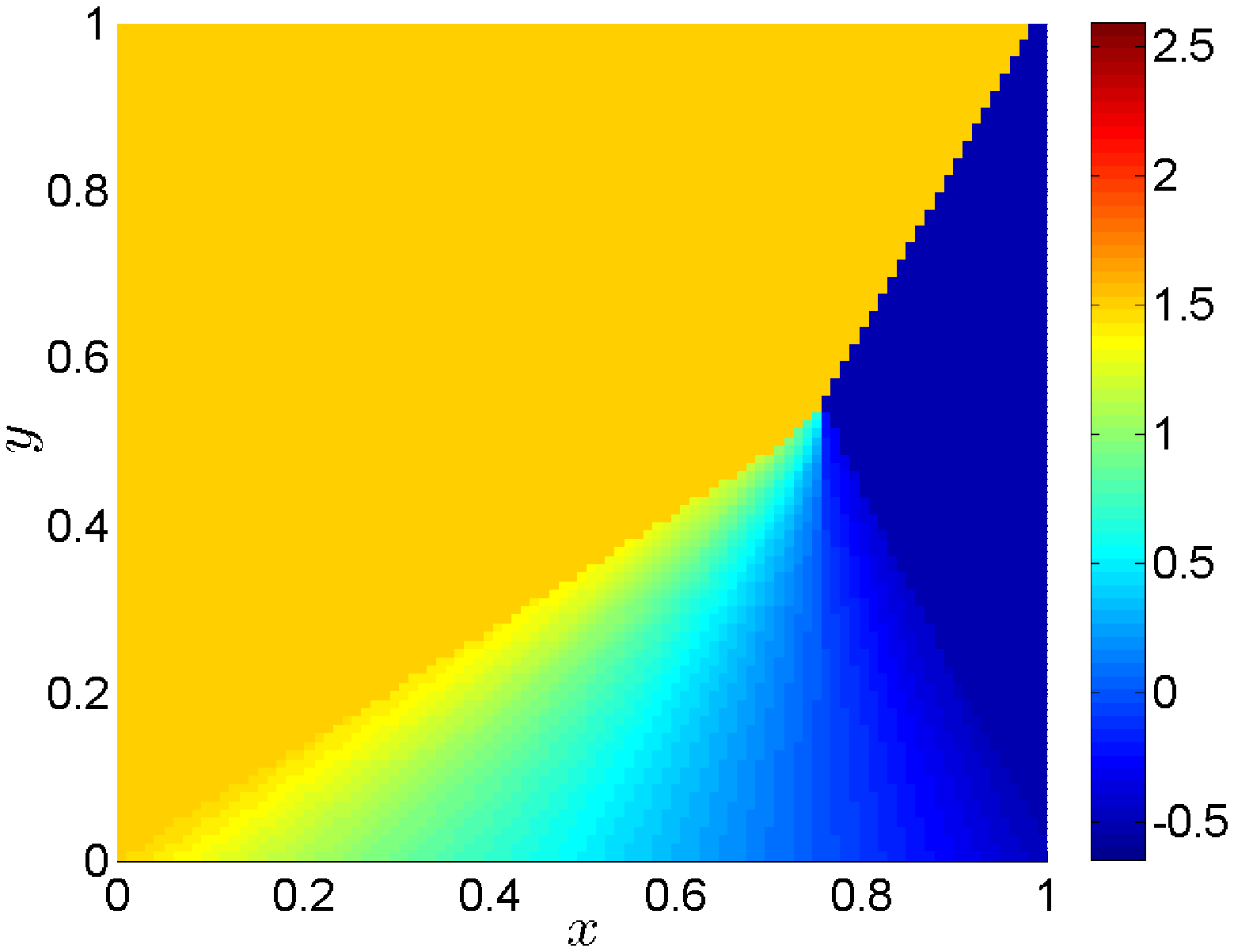}\label{fig:stupidLFSweep}}
\subfigure[]{\includegraphics[width=0.45\textwidth]{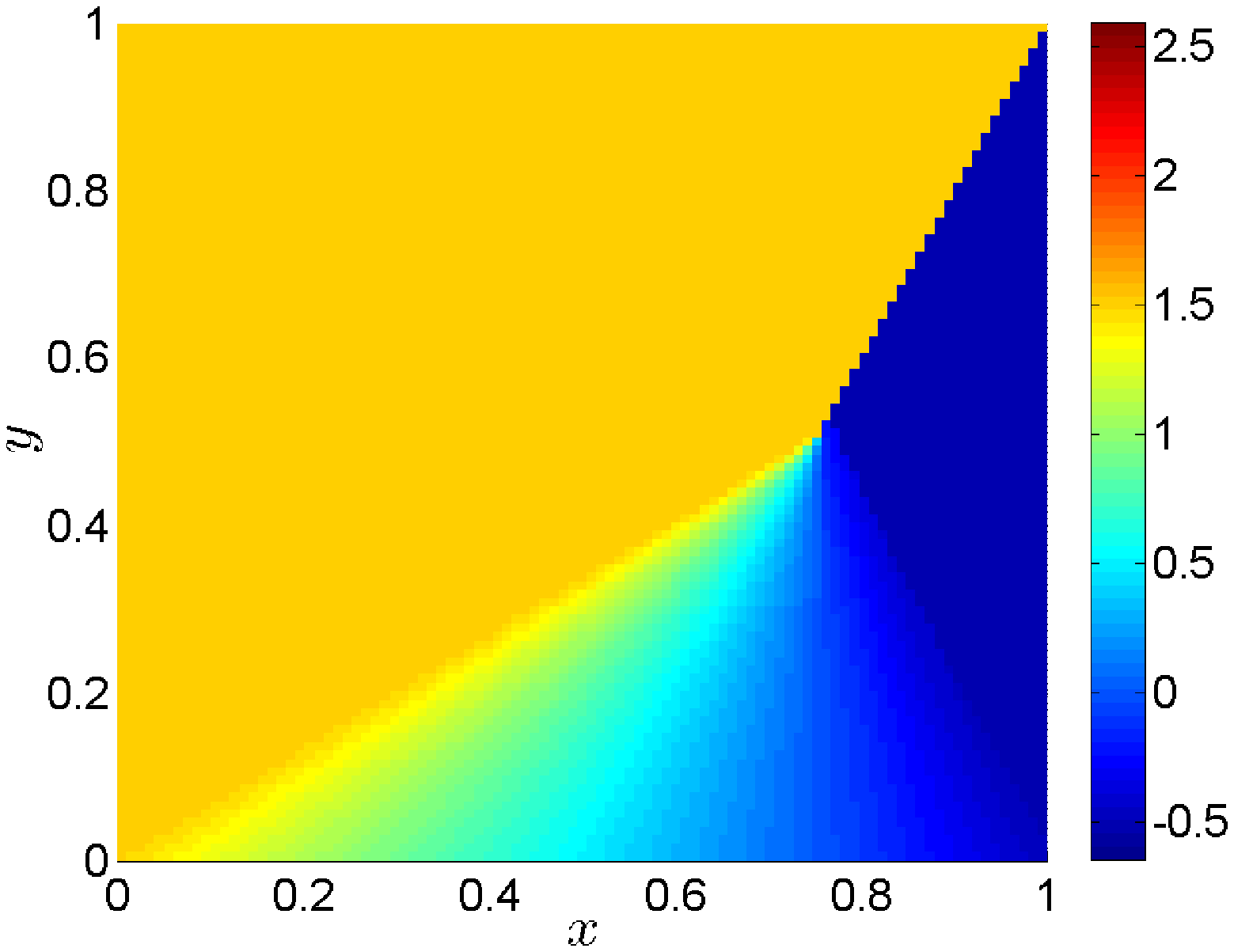}\label{fig:stupidLWSweep}}
\caption{Solution of~\eqref{eq:2dscalar} computed by shock-capturing using \subref{fig:stupidLF}~Lax-Friedrichs or \subref{fig:stupidLW}~Lax-Wendroff.  Solution computed by sweeping using \subref{fig:stupidLFSweep}~Lax-Friedrichs or \subref{fig:stupidLWSweep}~Lax-Wendroff.  All solutions are presented using the same colour axes, which are adjusted to allow for the overshoots in the Lax-Wendroff scheme.}
\label{fig:stupidCompare}
\end{center}
\end{figure}

The basic sweeping and matching technique involves growing a shock curve that begins at a boundary point, whereas in this example the shock begins in the interior of the domain.  Alternatively, we can simply sweep from the bottom using a conservative scheme.  However, this will produce a smeared-out shock (Figure~\ref{fig:stupidLF}).

Now we describe an alternative procedure, where the sweeping step is only used to generate smooth solution branches.  In particular, this means that we do not need to use a conservative scheme.

\begin{figure}[htdp]
\begin{center}
\subfigure[]{\includegraphics[width=0.45\textwidth]{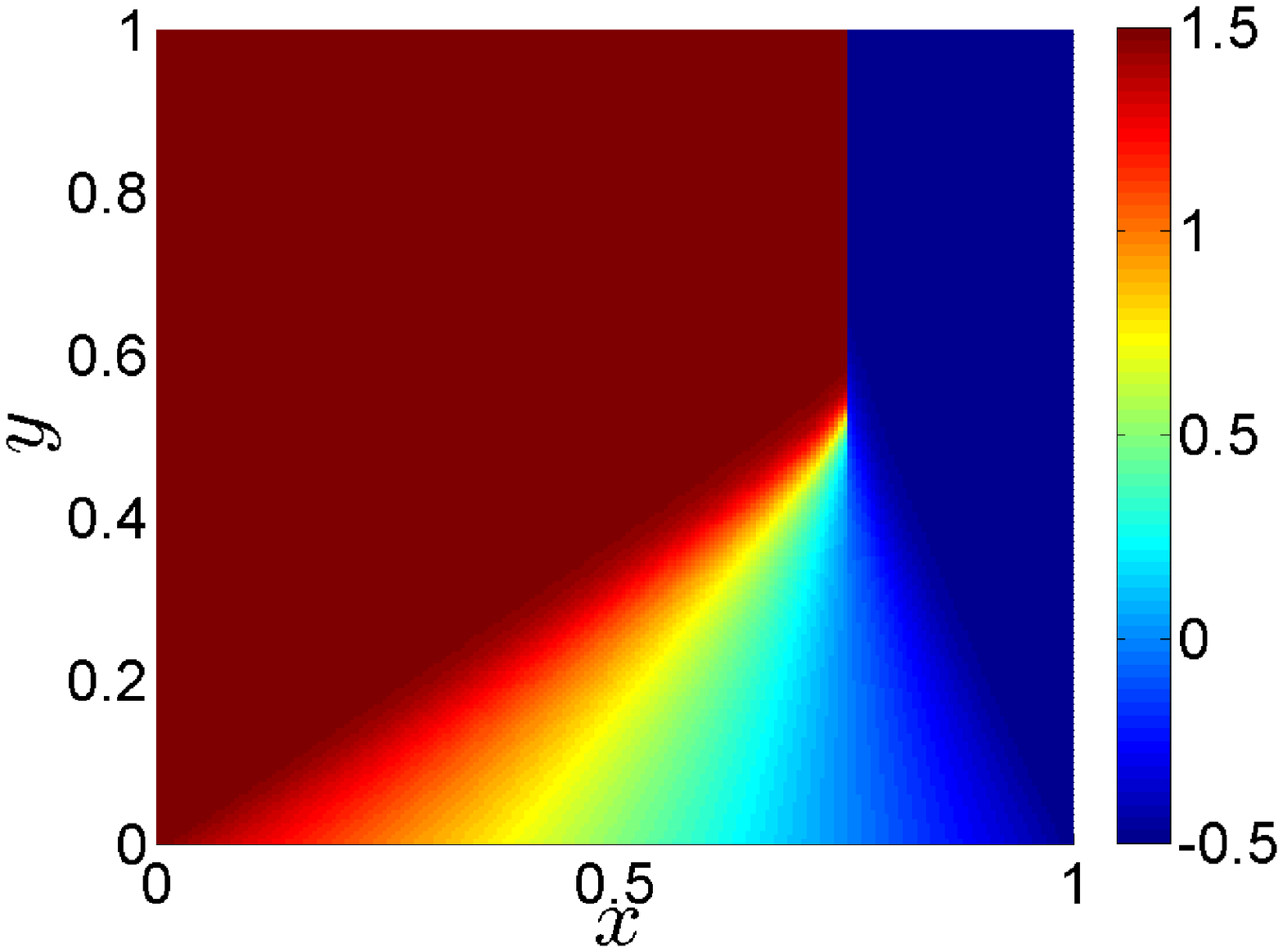}\label{fig:match_uB}}
\subfigure[]{\includegraphics[width=0.43\textwidth]{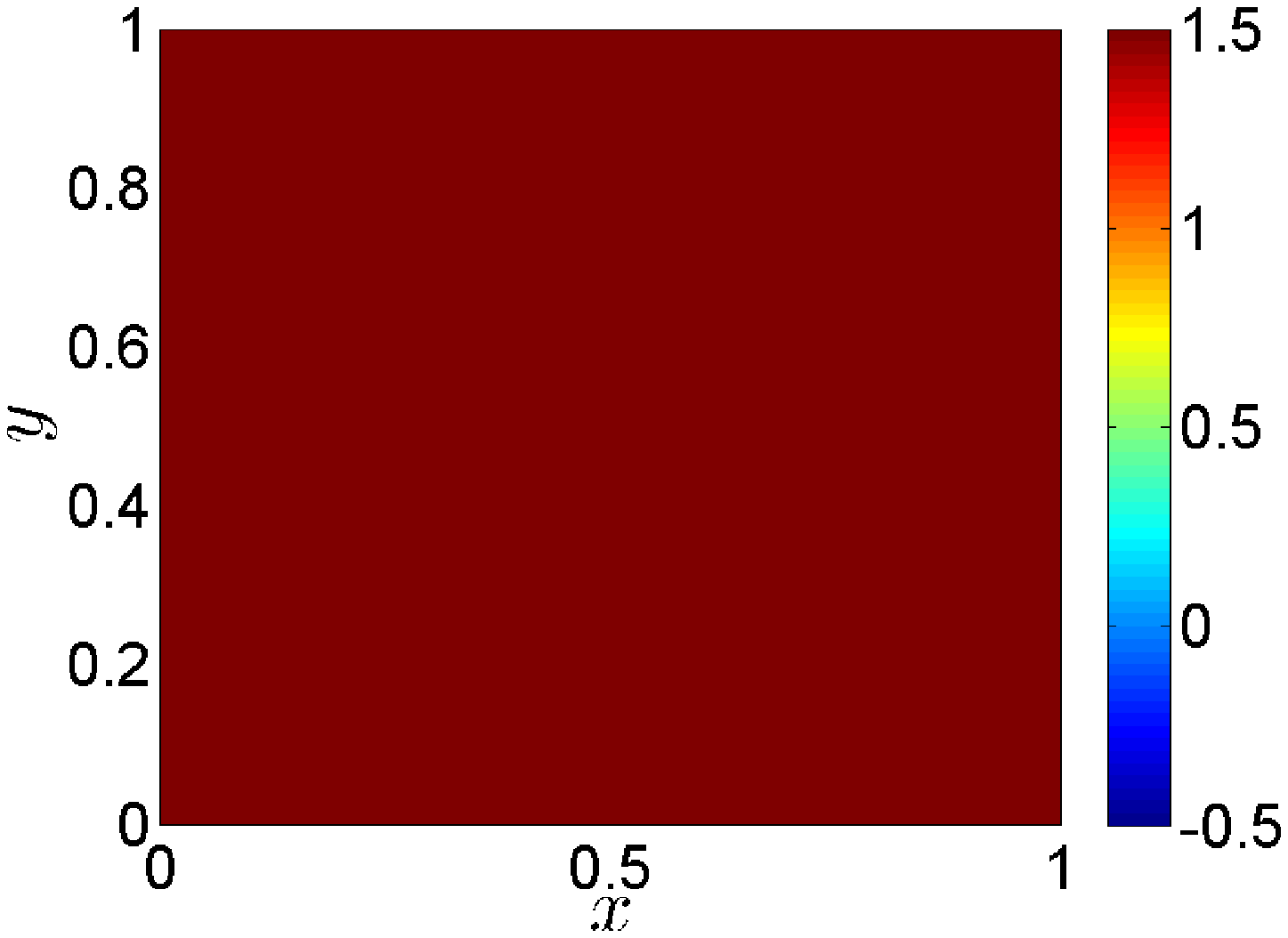}\label{fig:match_uL}}\\
\subfigure[]{\includegraphics[width=0.43\textwidth,height=0.34\textwidth]{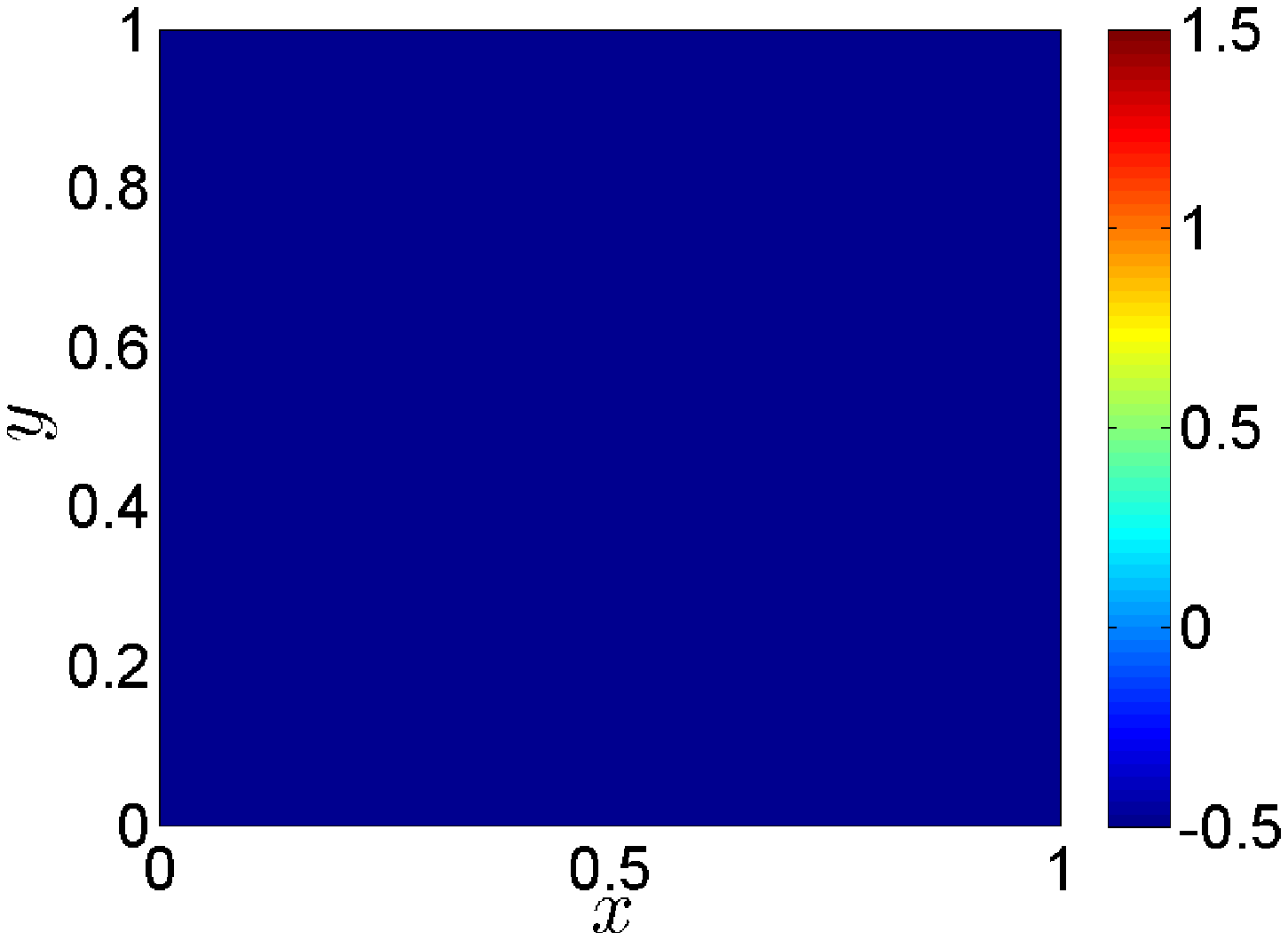}\label{fig:match_uR}}
\subfigure[]{\includegraphics[width=0.45\textwidth]{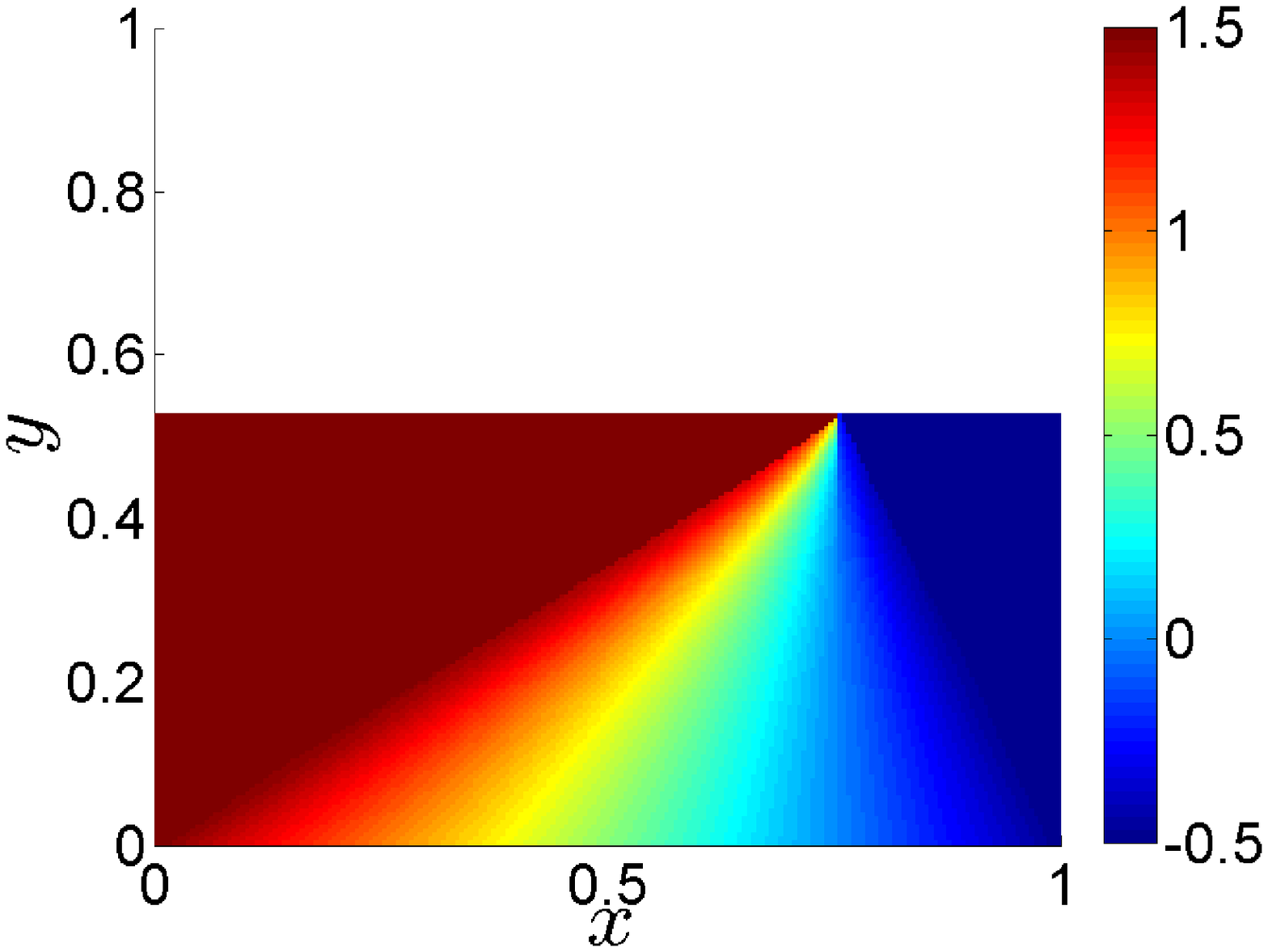}\label{fig:stupid_bottom}}\\
\subfigure[]{\includegraphics[width=0.45\textwidth]{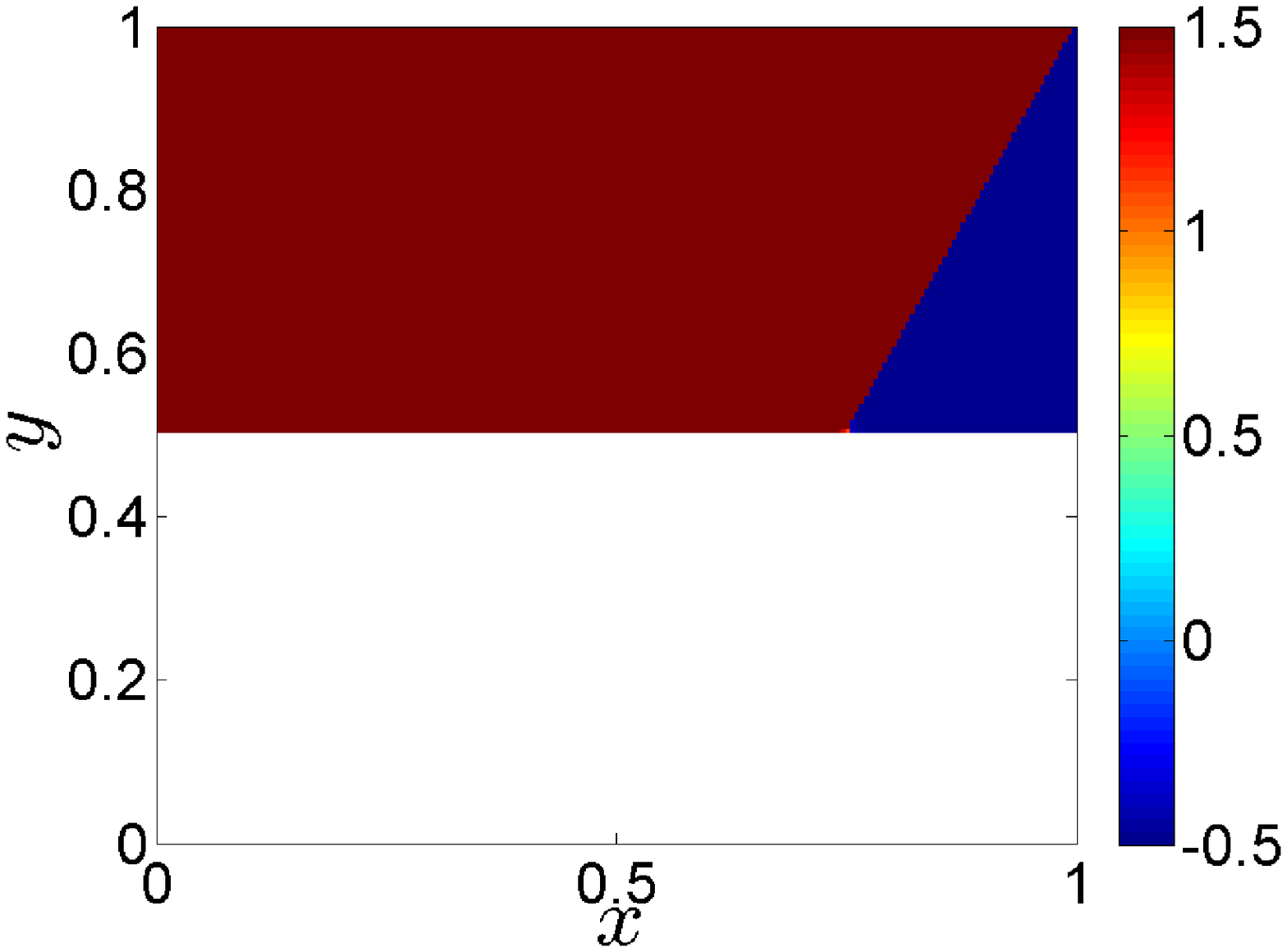}\label{fig:stupid_matchtop}}
\caption{Solution of~\eqref{eq:2dscalar}.  \subref{fig:match_uB}~Bottom, \subref{fig:match_uL}~left, and \subref{fig:match_uR}~right solution branches.  Matching done in the \subref{fig:stupid_bottom}~bottom and \subref{fig:stupid_matchtop}~top parts of the domain. 
}
\label{fig:match}
\end{center}
\end{figure}

First we need to solve the conservation law in the region below the shock.
\begin{enumerate}
\item We begin by identifying the point $x_*$ along the bottom boundary where the characteristics become vertical; this is the point where $f'(u(x^*,0)) = 0$.
\item We generate a bottom solution branch $u_B$.  This can be accomplished using any scheme, including a non-conservative scheme, since we will only make use of a smooth portion of the solution.  See Figure~\ref{fig:match_uB}.  Notice that we have chosen to use a non-conservative scheme, so that the resulting shock curve does not satisfy the Rankine-Hugoniot conditions.
\item Next we generate left and right solution branches $u_L$ and $u_R$ by sweeping from the sides.  
See Figure~\ref{fig:match_uL}-\ref{fig:match_uR}.
\item We compute a new solution branch by matching $u_B$ and $u_L$ in the region $x<x_*$.  This matching will start at the corner $(0,0)$ and continue until $x=x_*$, where the characteristics from the left boundary collide with the vertical characteristics from $(x_*,0)$.  We use $(x_*,y_*)$ to denote the endpoint of this computed curve.  
We similarly match this result with $u_R$ in the region $x>x*, y<y_*$.  See Figure~\ref{fig:stupid_bottom}.
\end{enumerate}
The result of this procedure is a solution $\tilde{u}$ that satisfies the conservation law in the region $y<y_*$ below the shock.

In order to compute the top portion of the solution, we need only solve the same PDE in the domain $[0,1]\times[y_*,1]$ subject to the boundary conditions
\[ u(0,y) = 1.5, \quad u(1,y) = -0.5, \quad u(x,y_*) = \tilde{u}(x,y_*).  \]
This problem now involves a shock originating from the boundary (at $(x_*,y_*)$), which can be solved using a straightforward application of the sweeping and matching techniques.  See Figure~\ref{fig:stupid_matchtop}.

Solutions constructed by this procedure are shown in Figure~\ref{fig:stupidCompare}.  To illustrate the advantages of the sweeping method, we use Lax-Friedrichs or Lax-Wendroff schemes to perform the sweeping steps.  For comparison, we also present solutions computing via shock-capturing with Lax-Friedrichs or Lax-Wendroff.  As desired, the sweeping method produces a clean shock curve, with none of the smearing or oscillations resulting from the capturing methods.

Computation times are displayed in Table~\ref{table:stupid} and validate the claim that the method has linear computational cost.

This example demonstrates several advantages of the fast sweeping approach.
\begin{itemize}
\item Sharp shock capturing: this result contains none of the smearing or oscillations that accompany a traditional shock-capturing method.
\item Flexibility: any consistent, stable scheme can be used to compute smooth solution branches.  In this example, a simple non-conservative scheme was used, but this did not affect the correctness of the computed shock curve.
\item Challenging shock structures can be handled by breaking the problem down into simpler problems in appropriate sub-domains.
\item The method is very computationally efficient.
\end{itemize}


\begin{table}[htdp]
\caption{Computation times on an $N\times N$ grid for the solution of~\eqref{eq:2dscalar} by sweeping.}
\begin{center}
\begin{tabular}{c||ccccc}
N  & 32& 64 & 128 & 256 & 512  \\
\hline
CPU Time (s) & 0.19 & 0.56 & 1.99 & 7.40 & 30.03 \\
\end{tabular}
\end{center}
\label{table:stupid}
\end{table}

\subsection{Rarefaction}\label{sec:rarefaction}
A second example where sweeping is desirable is in problems that involve rarefaction.
We consider the same equation as in the previous section,
\bq\label{eq:rarefaction} \left(\frac{u^2}{2}\right)_x + u_y = 0, \eq
in the domain $[-1,1]\times[0,1]$.
We enforce the boundary condition
\[ u(x,0) = \begin{cases} -1, & x<0 \\ 1/2, & x>0. \end{cases} \]
In this setting, the discontinuity at the boundary point $(x_*,0)=(0,0)$ will result in a rarefaction wave rather than a shock.

\begin{figure}[htdp]
\begin{center}
\subfigure[]{\includegraphics[width=0.45\textwidth]{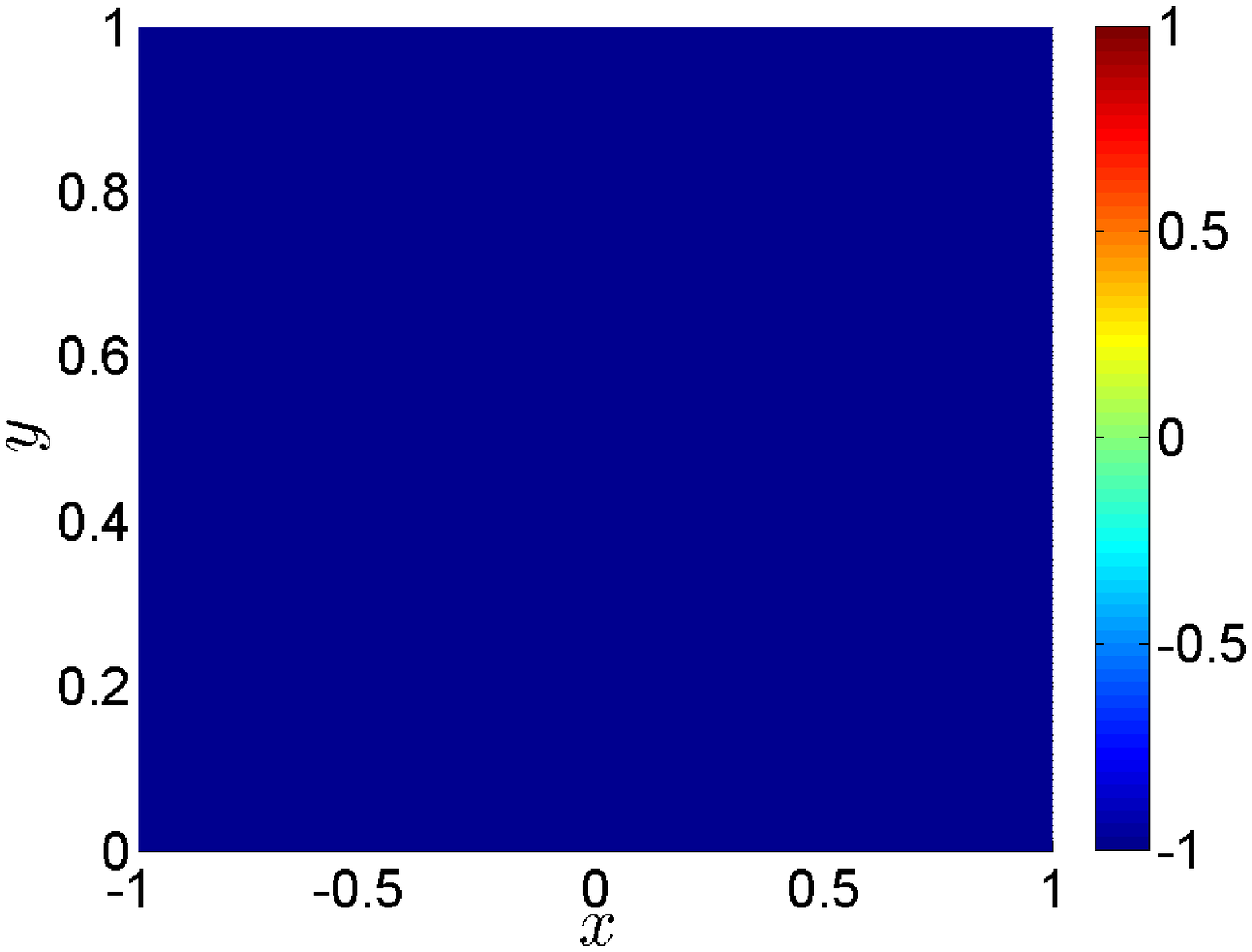}\label{fig:rarefactionUL}}
\subfigure[]{\includegraphics[width=0.45\textwidth]{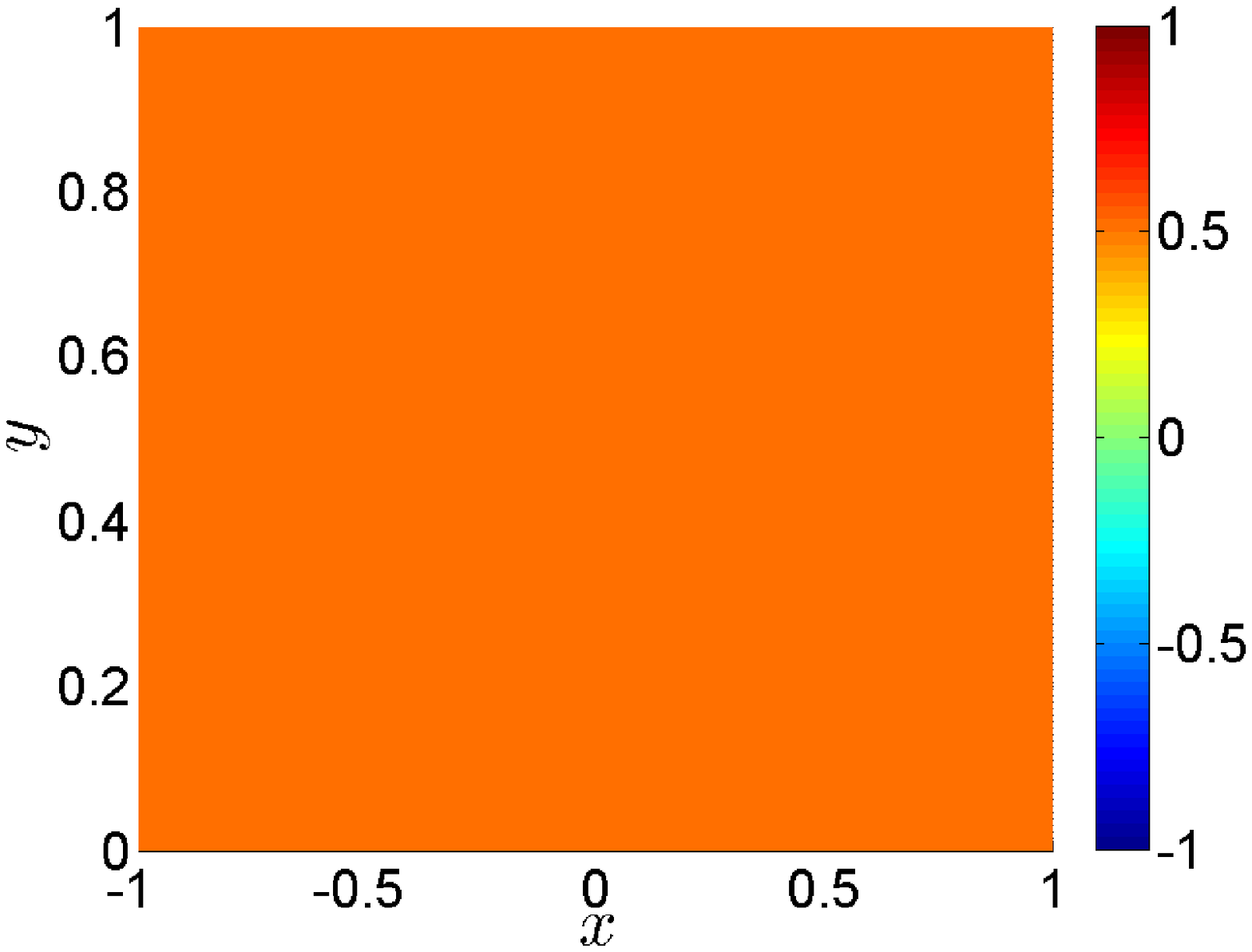}\label{fig:rarefactionUR}}\\
\subfigure[]{\includegraphics[width=0.45\textwidth]{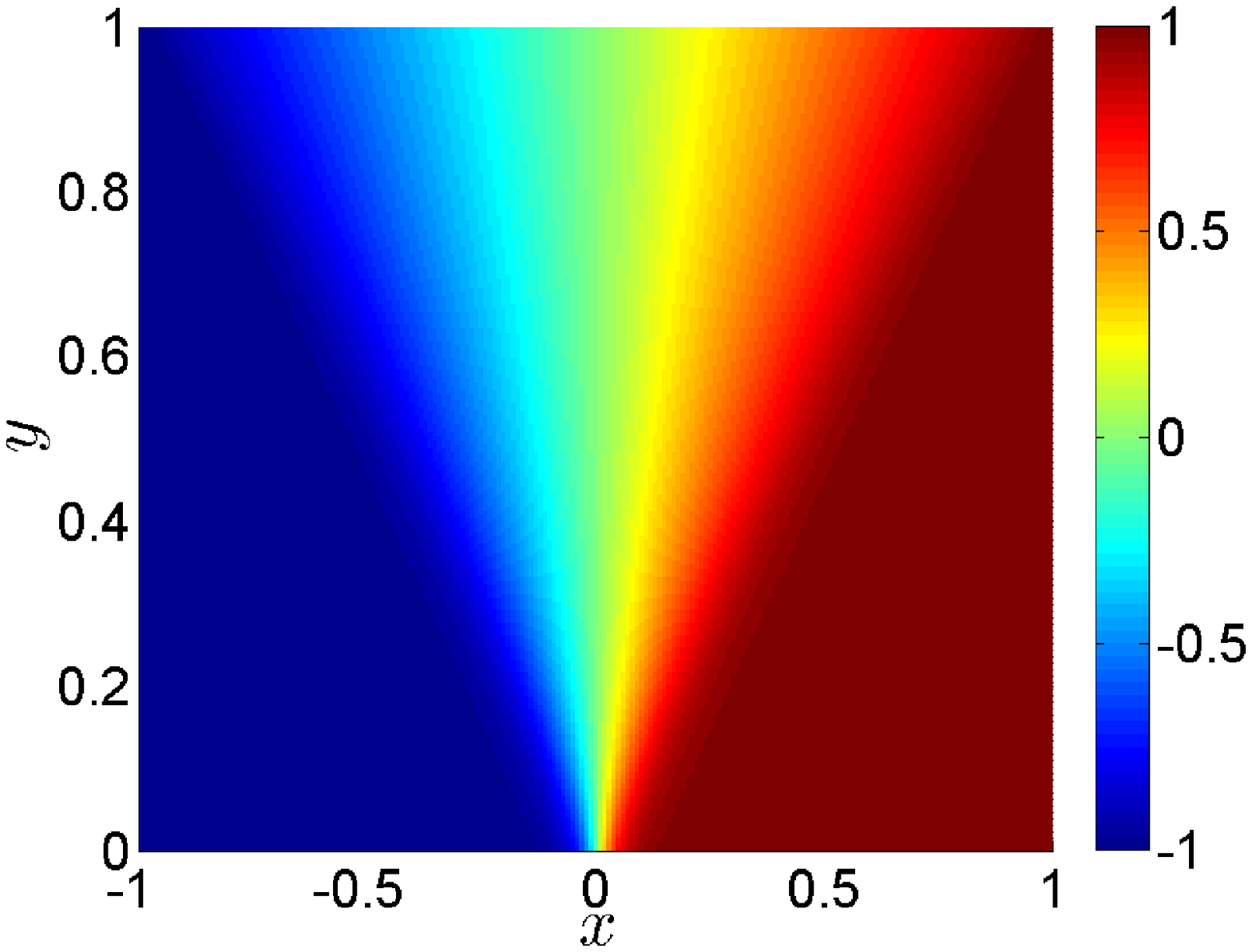}\label{fig:rarefactionUT}}
\subfigure[]{\includegraphics[width=0.45\textwidth]{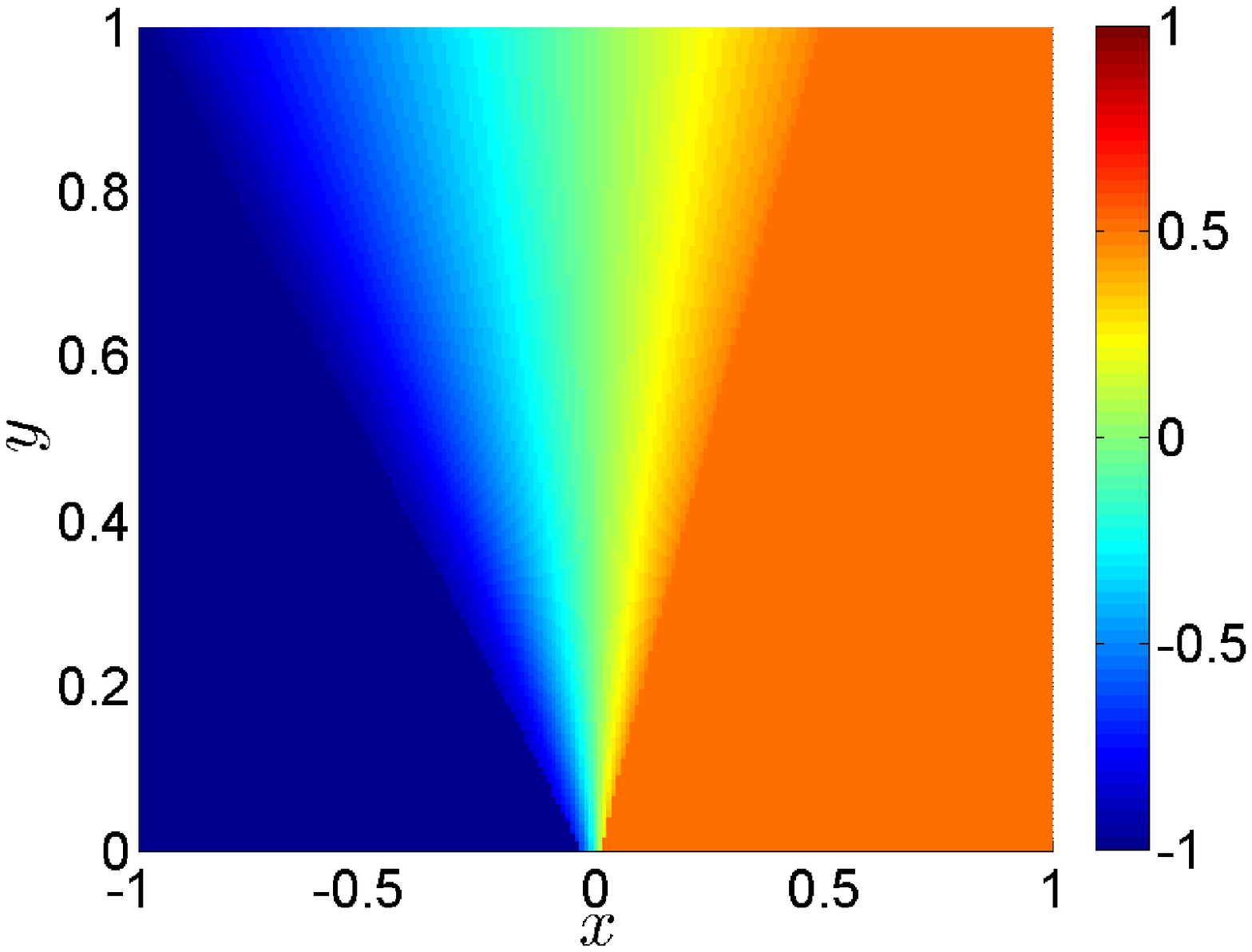}\label{fig:rarefactionU}}\\
\subfigure[]{\includegraphics[width=0.45\textwidth]{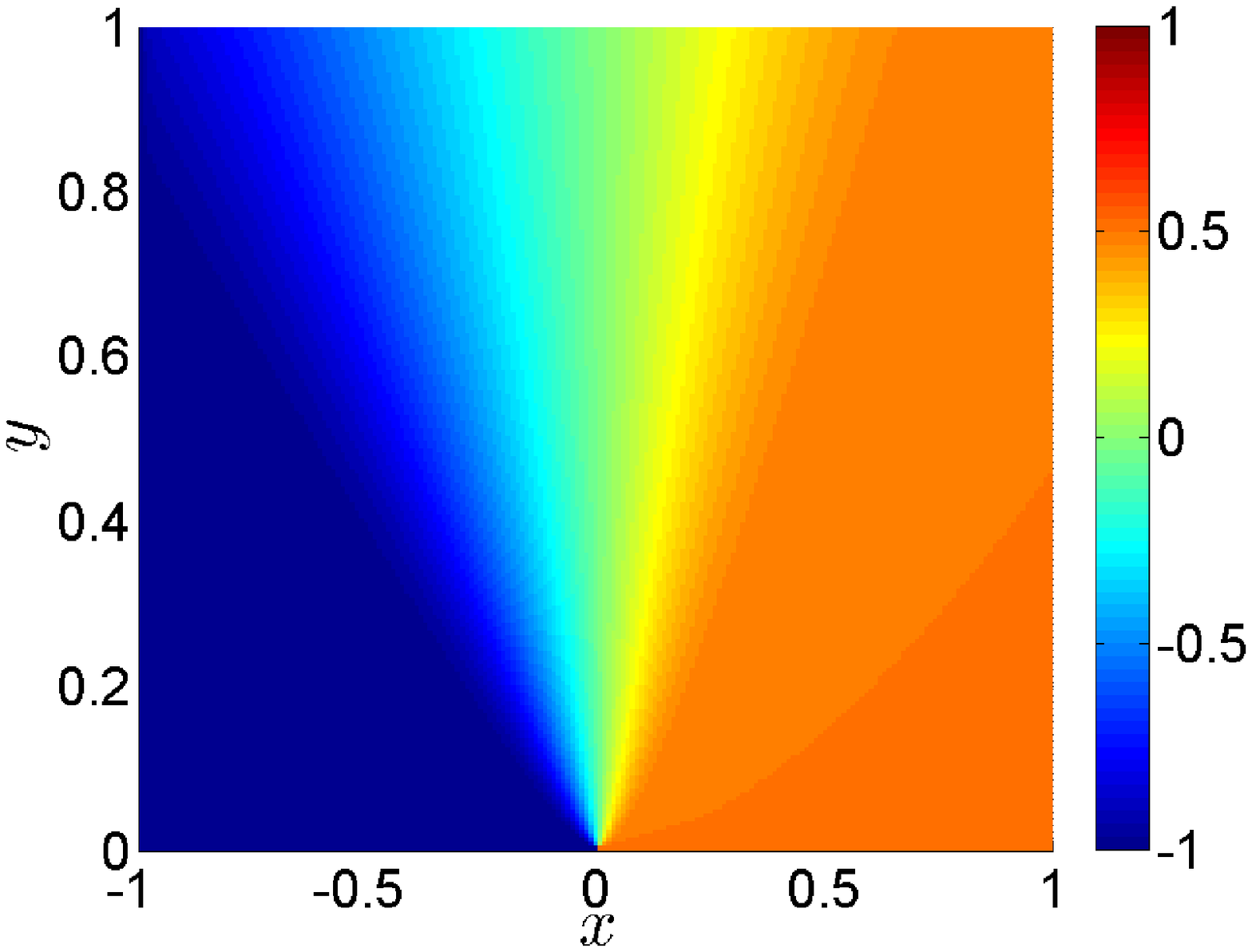}\label{fig:rarefactionULF}}
\caption{Solution of~\eqref{eq:rarefaction}.  \subref{fig:rarefactionUL}~Left, \subref{fig:rarefactionUR}~right, and \subref{fig:rarefactionUT}~top solution branches.  Solution computed using \subref{fig:rarefactionU}~sweeping and matching or \subref{fig:rarefactionULF}~a shock capturing scheme.}
\label{fig:rarefaction}
\end{center}
\end{figure}

As in the previous example, we could solve this by a single sweep from the bottom boundary as in Figure~\ref{fig:rarefactionULF}.  However, while this example does not include a shock curve, the solution is not differentiable.  Capturing schemes do not produce the sharp edge that is desirable.

We employ the following sweeping procedure.
\begin{enumerate}
\item Compute a left branch by sweeping from the bottom, using the given boundary conditions at $x<0, y=0$, and extending these smoothly along the remainder of the bottom boundary.  Similarly, compute a right solution branch.  See Figures~\ref{fig:rarefactionUL}-\ref{fig:rarefactionUR}.
\item Notice that these two solution branches cannot be matched across an entropy satisfying shock curve.  Instead, they must be connected through a rarefaction wave.
\item Compute a top branch by sweeping from the top boundary, solving the PDE in the form
\[-f(u)_x - g(u)_y = 0.\]  
To determine the correct boundary values at $(x,1)$, we recall that the characteristic emanating from this point should intersect the bottom boundary at the point $(x_*,0)$ where the given boundary conditions are discontinuous.  Additionally, the slope of this characteristic line will be given by $g'(u)/f'(u)$.  Consequently, the correct boundary value at each point $(x,1)$ can be determined by solving
\[ f'(u) = g'(u)(x-x_*). \]
See Figure~\ref{fig:rarefactionUT}.
\item Match the top and left solution branches starting from the point $(x_L,1)$ on the top boundary where the two solution branches agree.  In this setting, the matching curve is constructed so that the resulting solution will be continuous.  Similarly, match this result with the right solution branch.
\end{enumerate}

See Figure~\ref{fig:rarefactionU} for a picture of the solution computed by this procedure using a Lax-Friedrichs approximation of the flux.  This result is much sharper than the solution in Figure~\ref{fig:rarefactionULF}, which was computed using Lax-Friedrichs capturing.  This technique is also computationally efficient, as evidenced in Table~\ref{table:rarefaction}.

\begin{remark}
The first step of this process required us to produce a smooth extension of the boundary conditions into the region $x>0$.  In this example, a constant extension was natural.  When the boundary data is not constant, many extensions are possible.  To prevent dramatic changes in the orientation of the characteristics, it is desirable to minimise the variation in the boundary data.  One possible extension of more general boundary data $u(x,0)$ into the region $x>0$ can be obtained by computing
\[ \tilde{u}(x,0) = \begin{cases}
u(x,0) & x < 0\\
u(\epsilon,0) & x > \epsilon > 0\\
u(0,0) + \frac{x}{\epsilon}(u(\epsilon,0)-u(0,0)) & 0 < x < \epsilon. 
\end{cases} \]
This result can be convolved with a smoothing kernel in order to obtain a smooth extension.
\end{remark}

\begin{table}[htdp]
\caption{Computation times on an $N\times N$ grid for the solution of~\eqref{eq:rarefaction} by sweeping.}
\begin{center}
\begin{tabular}{c||ccccc}
N  & 32& 64 & 128 & 256 & 512  \\
\hline
CPU Time (s) & 0.24 & 0.73 & 2.83 & 10.98 & 42.48 \\
\end{tabular}
\end{center}
\label{table:rarefaction}
\end{table}

\section{Two-Dimensional Systems}\label{sec:systems}
Constructing sweeping methods for two-dimensional systems involves additional challenges that are not present in scalar problems.  One main challenge is the issue of ``incomplete'' boundary conditions, which is also present in one-dimensional systems.

\subsection{Shock reflection}\label{sec:reflection}
A setting where this occurs is the following shock reflection problem, which involves solving the stationary Euler equations
\bq\label{eq:euler} 
\left(\begin{tabular}{c}$\rho u$\\$\rho u^2+p$\\$\rho u v$\\$u(E+p)$\end{tabular}\right)_x  + \left(\begin{tabular}{c}$\rho v$\\$\rho u v$\\$\rho v^2+p$\\$v(E+p)$\end{tabular}\right)_y = 0  
\eq
in the domain
\[ 0 \leq x \leq 4, 0 \leq y \leq 1\]
with $p = (\gamma-1)\left(E-\frac{1}{2}\rho(u^2+v^2)\right)$ and $\gamma = 1.4$.  

Following~\cite{Chen_LFSweeping,Shu_WENOHomotopy}, we enforce the boundary conditions
\[
(\rho, u, v, p) = 
\begin{cases}
(1.69997,2.61934,-0.50632, 1.528191) & y = 1\\
(1, 2.9, 0, 1/\gamma) & x = 0.
\end{cases}
\]
A reflection condition (i.e. $v=0$) is imposed at $y=0$ and no boundary conditions are given at $x=4$.

In~\cite[Theorem~3]{EFTSweeping}, we showed that this problem could be solved using a single left-to-right sweep by solving a paraxial form of the equations.  In that setting, the shock is resolved using a shock capturing scheme.  As in the discussion in section~\ref{sec:matching}, this will produce a smeared out shock.  See Figure~\ref{fig:eulerCapturing}.

The initial steps of a fast sweeping method can be carried out on this problem.
\begin{enumerate}
\item Generate a smooth left state by sweeping in the boundary condition at $x=0$.
\item Generate a smooth top state by sweeping in the boundary condition at $y=1$.
\item Match the left and top states by constructing a shock curve that starts at the point $(0,1)$ where the boundary conditions are discontinuous.
\end{enumerate}

\begin{figure}[htdp]
	\centering
	{\includegraphics[width=0.9\textwidth]{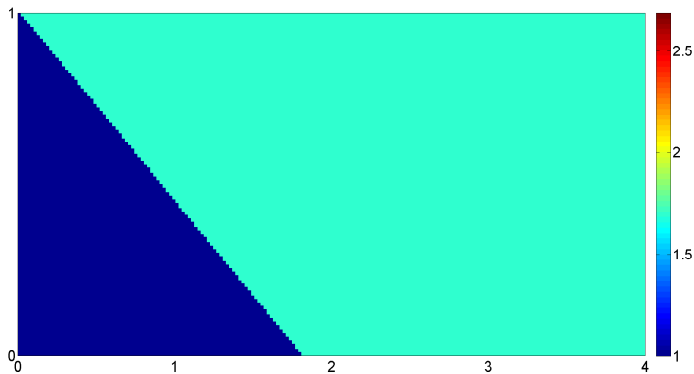}}
  	\caption{Density for the 2D Euler shock reflection problem  computed by matching top and left solution branches.}
  	\label{fig:rhoEulerTL}  	
\end{figure}

The result of this procedure is pictured in Figure~\ref{fig:rhoEulerTL}.  The shock separating the top and left branches eventually intersects the bottom boundary, as expected.  To the right of this intersection point $x_*$, the given bottom boundary condition ($v=0$) is not satisfied.  This indicates that another shock curve will need to extend from this boundary, separating the top state from a third right state.  However, determining this right state is non-trivial since only one boundary condition is given at the bottom, and none are given on the right side of the domain.

We do notice, though, that the initial right state $U_{right}$ and normal $n$ can be determined at the point $(x_*,0)$ where the shock originates.  This is because the reflection boundary condition gives us one component of the right state: $v_{right} = 0$. Thus the jump conditions
\[ \left(f(U_{top}(x))-f(U_{right}(x)), g(U_{top}(x))-g(U_{right}(x))\right) \cdot n = 0 \]
yield four equations for the normal $n$ and the three remaining components of $U_{right}$. 

The idea that we propose is to focus on the search for the unknown shock curve $y=\phi(x)$ that originates from $x_*$ and separates the top and right states.  Recall that we can easily compute the matched top state $U_{top}(x)$ throughout the domain.  If we are given this curve then we can employ the following procedure to construct the solution on the far side of the shock.
\begin{enumerate}
\item Apply the jump conditions 
\[ \left(f(U_{top}(x))-f(U_{right}(x)), g(U_{top}(x))-g(U_{right}(x))\right) \cdot n = 0 \]
at points on the shock curve, rejecting the continuous solution $U_{top}(x) = U_{right}(x)$.
\item Solve the system
\[
\begin{cases}
f(U)_x = - g(U)_y + a(U,x,y), & x>x_*, 0 < y < \phi(x)\\
U(x) = U_{right}(x), & x>x_*, y = \phi(x)\\
v(x) = 0, & x > x_*, y=0.
\end{cases}
\]
\end{enumerate}

\begin{figure}[htdp]
	\centering
			\subfigure[]{\includegraphics[width=.9\textwidth]{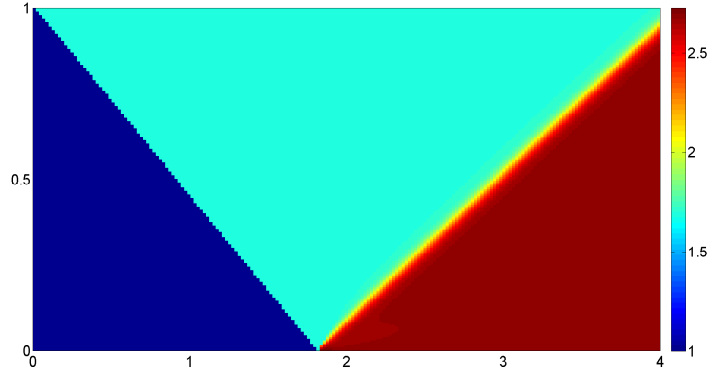}\label{fig:rho_Euler_wrongshock}}
			\subfigure[]{\includegraphics[width=.9\textwidth]{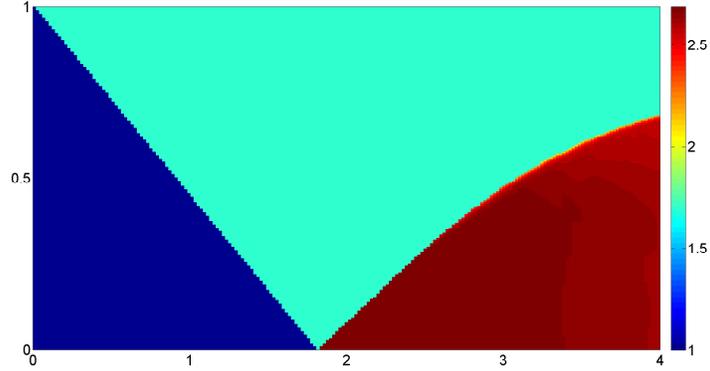}\label{fig:rho_Euler_curveshock}}
  	\caption{Density for the 2D Euler shock reflection problem computed using the sweeping approach with incorrect shock curves.}
  	\label{fig:EulerWrongShock}  	
\end{figure}

We can also apply the above procedure using an arbitrary shock curve $y=\phi(x)$.  However, in this case, the reflected solution coming from the boundary $y=0$ will be incompatible with the result of applying the jump conditions across the shock.  The result will be the formation of a sharp layer at some point in the domain, which contradicts the hypothesis that the data should produce a single reflected shock. See Figures~\ref{fig:rho_Euler_wrongshock}-\ref{fig:rho_Euler_curveshock} for examples of layers that form near two different candidate shock curves.  Thus regularity of the computed right solution branch is the condition that determines the correct shock curve.

Knowing the orientation of the shock curve near $(x_*,0)$, we turn our attention to extending this curve and right state $U_{right}$.
The approach for generating the shock curve and resulting right state is then to find the (discrete) curve $(x_i,\phi_i)$ such that a smoothness condition is satisfied,
\bq\label{eq:smoothness}
\Sm(\Phi(U_{top}(x_i);n_i)) =0, \quad x_* < x_i < x_R.
\eq
In the computations below, the smoothness indicator is a discrete estimate of the second derivative in the $y$-direction, which we set equal to zero.  For higher-order approximations, a higher-order derivative could be used.  This approach can be motivated by finite difference approximations, which are derived by truncating Taylor series, essentially assuming that higher derivatives are zero.    Alternatively, we could choose to minimise some indicator of solution smoothness.

In problems where all eigenvalues of $\nabla f$ are positive, it is not necessary to construct the entire curve $y=\phi(x)$ at once.  Instead, it can be generated one step at a time using a single left-to-right sweep.  Consequently, construction of the shock curve requires the solution of a sequence of nonlinear scalar equations for each $\phi_i$ rather then the solution of a nonlinear system for all values of $\phi$.

The approach is as follows.
\begin{enumerate}
\item To sweep the right solution branch and shock curve from $x_{i-1}$ to $x_{i}$ first guess at the new location of the shock curve $\phi_i$.
\item The normal to the curve can be estimated and the jump conditions applied to get the value of $U_{right}$ at $(x_i,\phi_i)$.
\item Do a single sweeping step to compute $U_{right}$ from $x_{i-1}$ to $x_i$ for all $0 < y_j < \phi_i$.  Enforce the reflection condition at the bottom and the correct values of $U_{right}$ at the shock.
\item If the value $(x_i,\phi_i)$ extends the shock curve with an incorrect orientation, it will not be possible to connect the top state $U_{top}$ to a single smooth right state $U_{right}$ via this shock curve.  Instead, the incorrect shock curve will trigger the formation of additional intermediate states, resulting in a breakdown in the regularity of $U_{right}$ that originates at the point $(x_{i-1},\phi_{i-1})$.  
\item Compute a smoothness indicator at $(x_i,\phi_i)$ and use this value to adjust the value of $\phi_i$ via Newton's method.  If an undesirable intermediate state has originated at $(x_{i-1},\phi_{i-1})$, discrete approximations of derivatives in the $y$-direction will be large near $(x_i,\phi_i)$.  Thus $u_{yy}$ emerges as a natural indicator of smoothness.    
\end{enumerate}

\begin{figure}[htdp]
	\centering
	\subfigure[]{\includegraphics[width=0.9\textwidth]{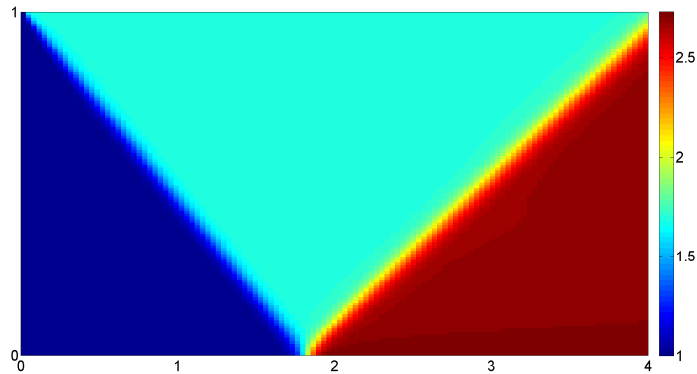}\label{fig:eulerCapturing}}
	\subfigure[]{\includegraphics[width=0.9\textwidth]{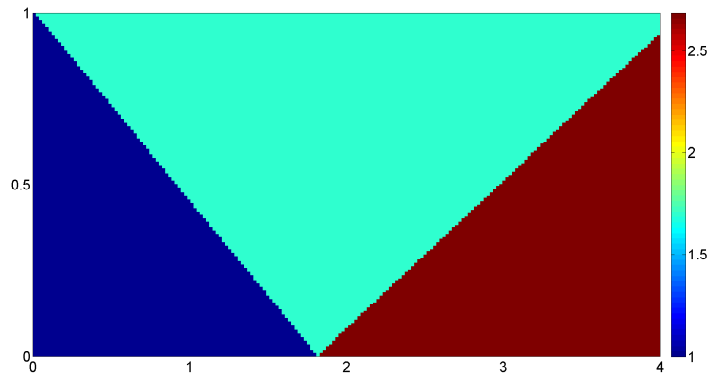}\label{fig:rho_Euler}}
  	\caption{Density for the 2D Euler shock reflection problem  computed using \subref{fig:eulerCapturing}~a shock-capturing method to solve a paraxial form of the equations and \subref{fig:rho_Euler}~sweeping.}
  	\label{fig:eulerCapturing}  	
\end{figure}

The solution computed using this fast sweeping approach is pictured in Figure~\ref{fig:rho_Euler}.  It contains a sharp shock curve, as desired.  We emphasise again that computing this solution was done efficiently, requiring only three sweeps through the data. Computation times are presented in Table~\ref{table:reflection} and validate our claims of optimal computational complexity.

\begin{table}[htdp]
\caption{Computation times on an $N\times N$ grid for the shock reflection problem.}
\begin{center}
\begin{tabular}{c||ccccc}
N  & 32& 64 & 128 & 256 & 512  \\
\hline
CPU Time (s) & 6.4 & 24.4 & 85.7 & 301.6 & 1096.8 \\
\end{tabular}
\end{center}
\label{table:reflection}
\end{table}

\subsection{Oblique shock (piecewise constant)}\label{sec:oblique}

We consider a second example of a two-dimensional system and use the same procedure described in the previous section.  In this oblique shock problem, a uniform flow impinges on an impermeable wall.  Following~\cite{Oblique}, we pose the two-dimensional stationary Euler equations in the domain.
\[\{(x,y)\in[0,1]\times[0,1/2]\mid y \geq (x-1/2)\tan\delta\}. \]
A horizontal flow is enforced at the left boundary,
\[ (\rho, u, v, p) = (1,3,0,1/\gamma), \quad x = 0 \]
and the bottom boundary is an impermeable wall,
\[ v = \begin{cases}
0, & x < 1/2, y = 0\\
u\tan\delta, & x > 1/2, y = (x-1/2)\tan\delta.
\end{cases} \]
When the uniform horizontal flow impinges on the wedge at the bottom of the domain, a shock will form.  

We solve this problem setting $\delta=15^\circ$  as the angle the bottom wedge makes with the horizontal.  We expect the shock to make an angle of approximately $32.2^\circ$ with the horizontal, and the flow beyond the shock should have a Mach number of approximately 2.255~\cite{Oblique}.

Solution of this problem via fast sweeping requires a two-step process:
\begin{enumerate}
\item Generate a smooth left state by sweeping in the boundary conditions at $x=0$.
\item Use the method described in \autoref{sec:reflection} to construct a shock curve and smooth right state emanating from the point $(1/2,0)$ where the bottom boundary condition is discontinuous.
\end{enumerate}

\begin{figure}[htdp]
  \centering
	\includegraphics[width=0.9\textwidth]{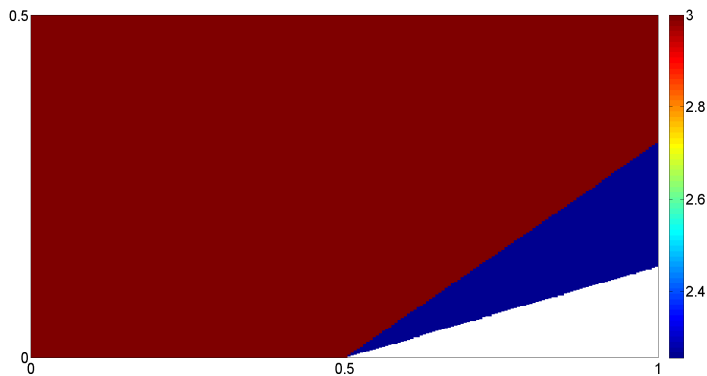}
	\caption{Mach number computed using sweeping for the oblique shock problem.}
	\label{fig:mach}
\end{figure}

The computed Mach number is displayed in Figure~\ref{fig:mach}.  In particular, the right state has a Mach number of approximately 2.2549, while the computed shock makes an angle of $32.240^\circ$ with the horizontal.  This is in agreement with the expected solution.  We also observe optimal computational complexity; see Table~\ref{table:oblique}.

\begin{table}[htdp]
\caption{Computation times on an $N\times N$ grid for the oblique shock problem with piecewise constant data.}
\begin{center}
\begin{tabular}{c||ccccc}
N  & 32& 64 & 128 & 256 & 512  \\
\hline
CPU Time (s) & 8.1 & 30.1 & 109.2 & 331.9 & 1334.3 \\
\end{tabular}
\end{center}
\label{table:oblique}
\end{table}

\subsection{Oblique Shock (nonconstant)}
Because the examples presented in the previous subsection involved piecewise constant data, the resulting shock curves were linear.  To better demonstrate the abilities of the fast sweeping method, we repeat the oblique shock wave problem, this time imposing non-constant data at the left boundary.  The flow at $x=0$ is now given by
\[ (\rho, u, v, p) = \left(1,3,10y\left(\frac{1}{2}-y\right),1/\gamma-0.3\sin(4\pi y)\right), \quad x=0. \]

\begin{figure}[htdp]
	\centering
       \subfigure[]{\includegraphics[width=.9\textwidth]{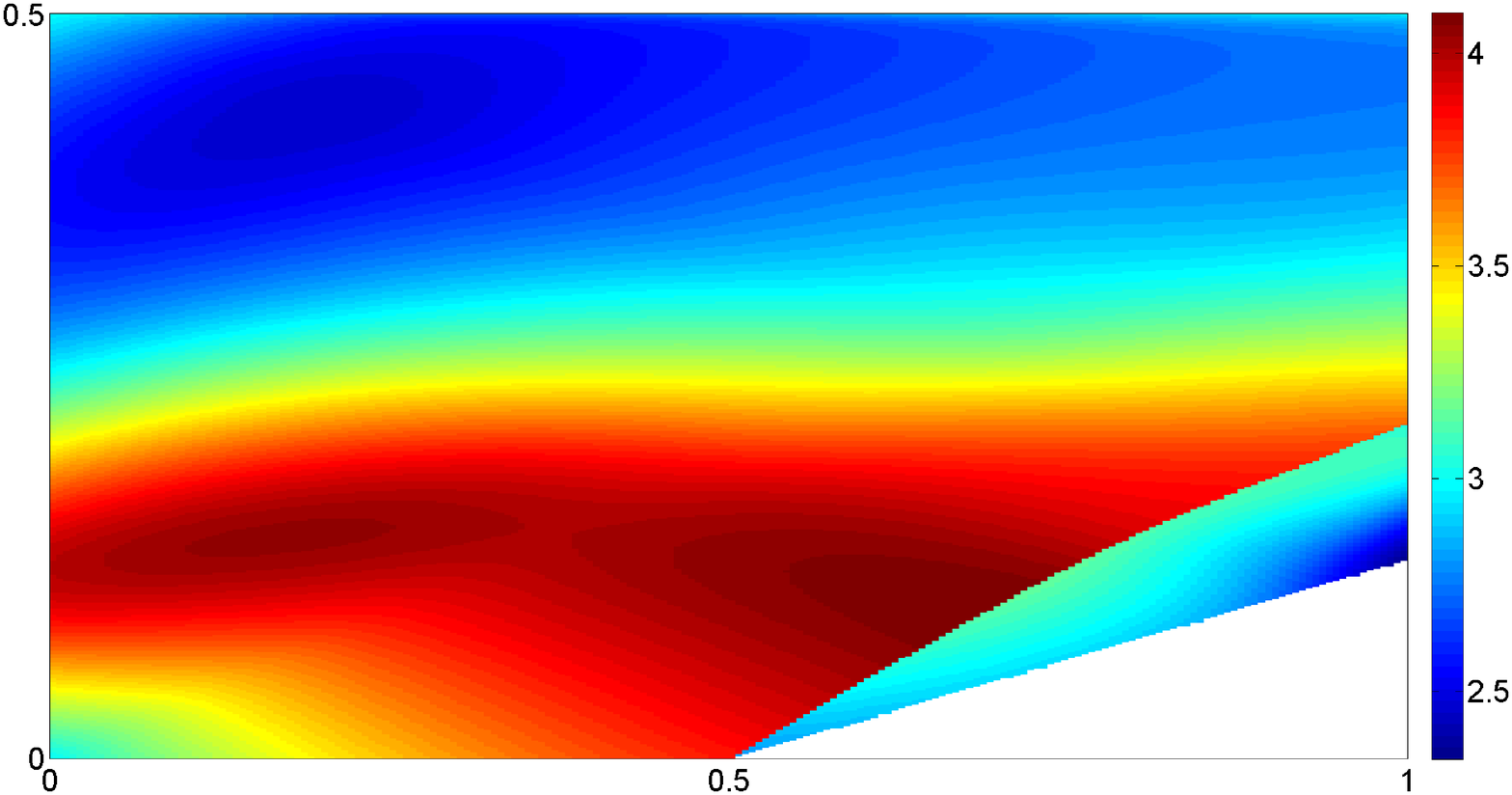}\label{fig:oblique_nonconstant2}}\\
			\vspace*{-12pt}
			\subfigure[]{\includegraphics[width=.9\textwidth]{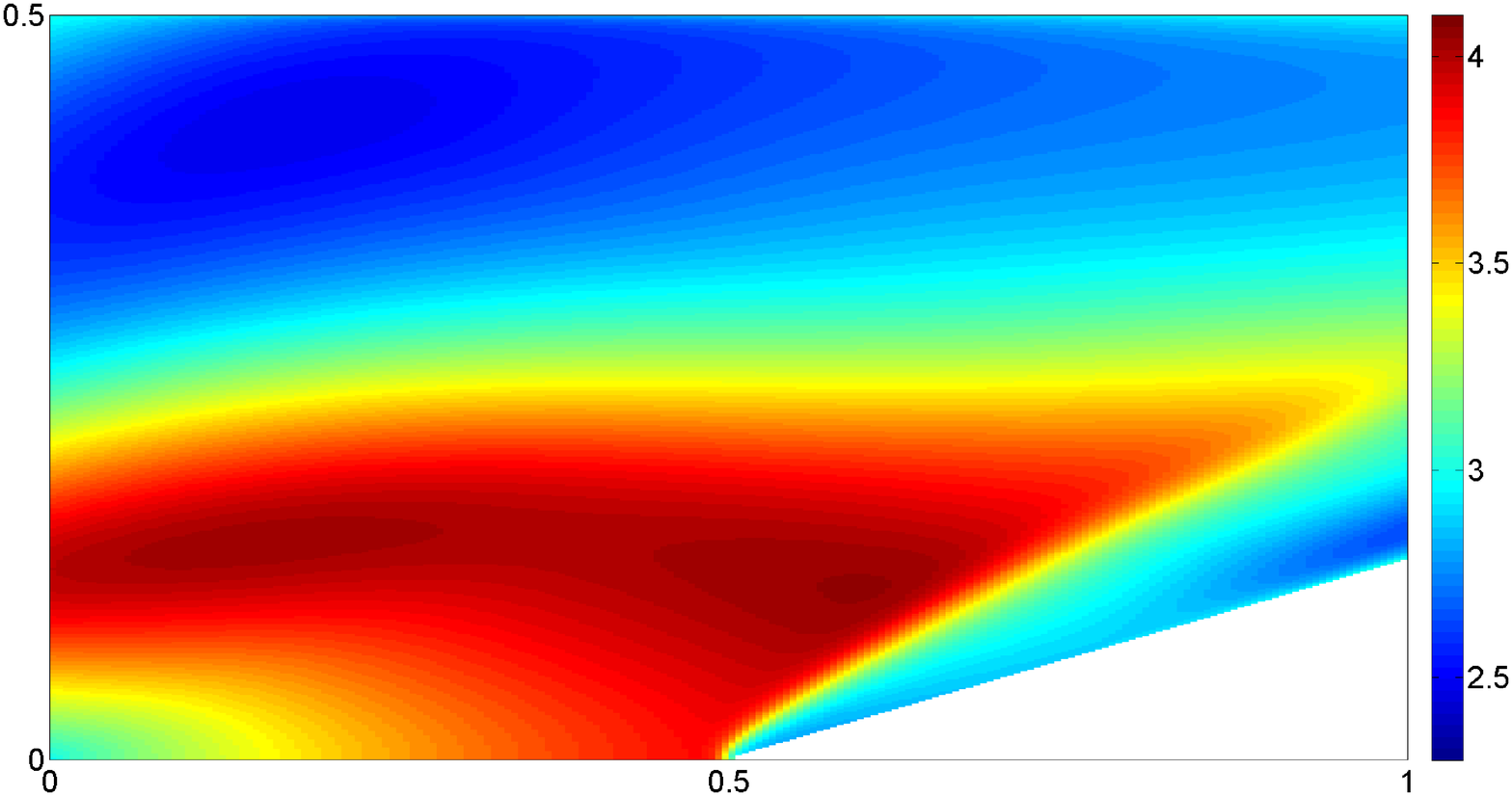}\label{fig:oblique_nonconstant_evolveC}}\\
			\vspace*{-12pt}
			\subfigure[]{\includegraphics[width=.9\textwidth]{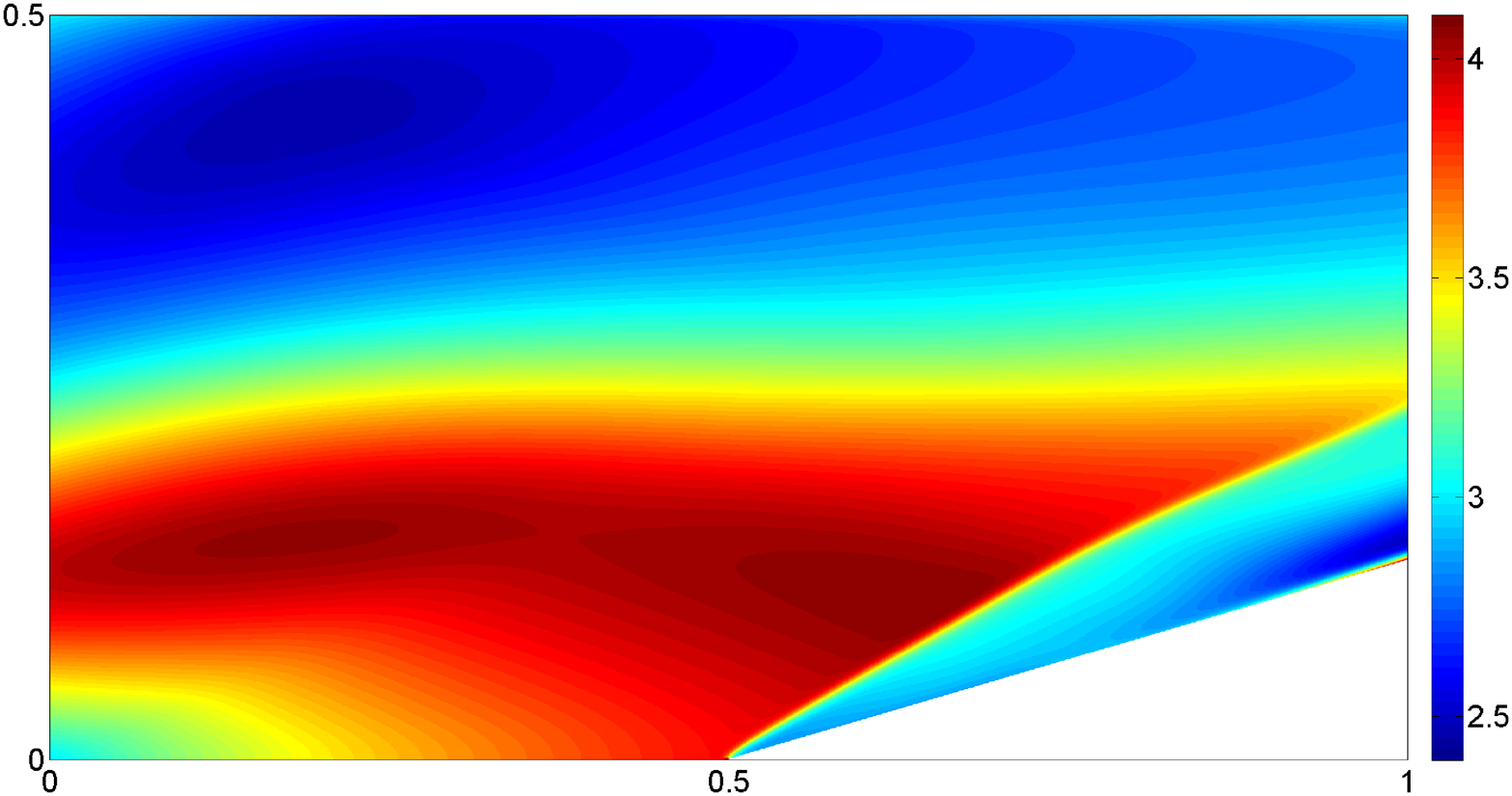}\label{fig:oblique_nonconstant_evolve}}
  	\caption{Mach number for the oblique shock problem with nonconstant data.  Solution obtained by~\subref{fig:oblique_nonconstant2}~sweeping or \subref{fig:oblique_nonconstant_evolveC}~evolving a Lax-Friedrichs scheme to steady state on a $200\times200$ grid. \subref{fig:oblique_nonconstant_evolve}~Reference solution computed by evolving a Lax-Friedrichs scheme to steady state on a $801\times801$ grid.}  	
  	\label{fig:oblique_nonconstant}
\end{figure}

The computed Mach number, which now includes variable states and a nonlinear shock curve, is displayed in Figure~\ref{fig:oblique_nonconstant2}.  As in the previous example, this solution was computed using a Lax-Friedrichs approximation, but produces much sharper results than a solution computed by simply evolving a Lax-Friedrichs scheme to steady state (Figure~\ref{fig:oblique_nonconstant_evolveC}).  The computational cost of the sweeping method is no greater than it was for the simpler piecewise constant case (Table~\ref{table:obliqueNC}).

No exact solution is available for this problem.  For reference, we include a plot of the Mach number obtained by evolving a Lax-Friedrichs scheme to steady state on a more refined ($801\times801$) grid (Figure~\ref{fig:oblique_nonconstant_evolve}).

\begin{table}[htdp]
\caption{Computation times on an $N\times N$ grid for the oblique shock problem with non-constant data.}
\begin{center}
\begin{tabular}{c||ccccc}
N  & 32& 64 & 128 & 256 & 512  \\
\hline
CPU Time (s) & 7.5 & 23.8 & 81.5 & 299.0 & 1147.2 \\
\end{tabular}
\end{center}
\label{table:obliqueNC}
\end{table}

\section{Conclusions}\label{sec:conclusions}
In this article, we have further developed a fast sweeping approach for computing steady state solutions to systems of conservation laws, which was first introduced in~\cite{EFTSweeping}.  We demonstrated the possibility of constructing methods that achieve higher accuracy than is possible with traditional shock-capturing methods.  We also extended the use of two-dimensional sweeping methods to problems that have a more complicated shock structure, including shocks that begin in the interior of the domain and rarefaction.  Finally, we developed fast sweeping methods for a class of two-dimensional systems.  In all cases, our method produced results with correct, sharp shock curves.  Computational experiments also validated the claim that the methods have optimal computational complexity.

Future challenges include the extension of these techniques to higher dimensions.  In problems with more challenging structures, it may not be straightforward to generate an explicit formulation of the shock surface.  An alternative
approach would be to represent the shock implicitly; for example, using a level set formulation.

\section*{Acknowledgements}
Yen-Hsi Richard Tsai thanks the National Center for Theoretical Sciences, Taipei, Taiwan for hosting his visit to Taipei,
where the research for this article was initiated.  He also thanks I-Liang Chern for stimulating conversation on related topics.

\end{document}